\documentclass[a4paper,12pt,twoside]{article}
\usepackage[latin1]{inputenc}                 
\usepackage{amsfonts,amssymb,amsmath,exscale} 
\usepackage[dvips]{graphicx,color,psfrag}     
\usepackage{fancyhdr}                         
\usepackage{caption}                          
\usepackage{rotating}
\usepackage{hyperref}
\usepackage{defcml}                           
\usepackage[numbers]{natbib}                  
\Floatboxname{Box}
\usepackage{verbatim}                         

\usepackage{amsmath,amsthm,amsfonts,amssymb,amscd}
\usepackage{subfig}
\usepackage{graphicx,bm,color}
\usepackage{algorithm}
        \usepackage{setspace}
\usepackage{algpseudocode}
\usepackage{stmaryrd}
\usepackage[framed,numbered,autolinebreaks,useliterate]{mcode}
\usepackage{blindtext}
\usepackage{multicol}
\usepackage{xargs}                      
\usepackage[color=gray!20, backgroundcolor=yellow, textwidth=1.85cm]{todonotes}
\newcommandx{\ju}[2][1=]{\todo[linecolor=blue,backgroundcolor=blue!25,bordercolor=blue,#1]{#2}}

\usepackage{mathtools}

\usepackage{booktabs,ctable,multirow,longtable} 
\usepackage[latin1]{inputenc}                 
\usepackage{amsfonts,amssymb,amsmath,exscale} 
\usepackage[dvips]{graphicx,color,psfrag}     
\usepackage{fancyhdr}                         
\usepackage{caption}                          
\usepackage{rotating}
\usepackage{defcml}                            
\usepackage[numbers]{natbib}                  
\usepackage{url}
\Floatboxname{Box}
\newcolumntype{\xx}{\bm{x}}
\renewcommand{\P}{\mathbb{P}}
\newcommand{\xx}{\boldsymbol{x}}
\newcommand{\R}{\mathbb{R}}

\DeclareMathOperator{\RMSE}{RMSE}
\DeclareMathOperator{\EMC}{E_{MC}}
\newcommand{\E}{\mathbb{E}}

\newtheorem{theorem}{Theorem}

\newtheorem{lemma}[theorem]{Lemma}
\newtheorem{proposition}{Proposition}

\definecolor{gray}{gray}{0.6}

\newtheorem{Remark}{Remark}[section]

\def\R{\mathbb R}

\usepackage{graphicx}

\fancypagestyle{plain}{%
	\fancyhead{}%
	\fancyfoot[c]{\sffamily\thepage}%
}
\makeatletter
 
\usepackage{scalerel,stackengine}
\stackMath
\newcommand\reallywidecheck[1]{%
	\savestack{\tmpbox}{\stretchto{%
			\scaleto{%
				\scalerel*[\widthof{\ensuremath{#1}}]{\kern-.6pt\bigwedge\kern-.6pt}%
				{\rule[-\textheight/2]{1ex}{\textheight}}
			}{\textheight}%
		}{0.5ex}}%
	\stackon[1pt]{#1}{\scalebox{-1}{\tmpbox}}%
}

\begin{document}

\thispagestyle{empty}
\vspace{-1cm}
\ce{\bf\Large Probabilistic failure mechanisms via Monte Carlo}
\vspace*{0.25cm}
\ce{\bf\Large  simulations of complex microstructures}

\vskip .35in

\ce{
Nima Noii\(^{a}\),	Amirreza Khodadadian\(^{b}\), Fadi Aldakheel\(^{a,c,}\)\footnote{Corresponding author.\\[1mm] 
E-mail addresses
khodadadian@ifam.uni-hannover.de(A. Khodadadian);	noii@ikm.uni-hannover.de (N. Noii);
aldakheel@ikm.uni-hannover.de (F. Aldakheel).
}}
\vskip .25in
\vspace{-0.2cm}
\ce{\(^a\) Institute of Continuum Mechanics} \ce{Leibniz Universit\"at Hannover, An der Universit\"at 1, 30823 Garbsen, Germany}\vskip .25in
\vspace{-0.2cm}
\ce{\(^b\) Institute of Applied Mathematics} \ce{Leibniz Universit\"at Hannover, Welfengarten 1, 30167 Hannover, Germany} \vskip .22in
\ce{\(^c\) Zienkiewicz Centre for Computational Engineering, Faculty of Science and Engineering} \ce{Swansea University, Bay Campus, SA1 8EN, UK} \vskip .22in

\begin{Abstract}
	
A probabilistic approach to phase-field brittle and ductile fracture with random material and geometric properties is proposed within this work. In the macroscopic failure mechanics, materials properties and exactness of spatial quantities (of different phases in the geometrical domain) are assumed to be homogeneous and deterministic. This is unlike the lower-scale with strong fluctuation in the material and geometrical properties. Such a response is approximated through some uncertainty in the model problem. The presented contribution is devoted to providing a mathematical framework for modeling uncertainty through stochastic analysis of a microstructure undergoing brittle/ductile failure. Hereby, the proposed model employs various representative volume elements with random distribution of stiff-inclusions and voids within the composite structure. We develop an allocating strategy to allocate the heterogeneities and generate the corresponding meshes in two- and three-dimensional cases. Then the Monte Carlo finite element technique is employed for solving the stochastic PDE-based model and approximate the expectation and the variance of the solution field of brittle/ductile failure by evaluating a large number of samples. For the prediction of failure mechanisms, we rely on the phase-field approach which is a widely adopted framework for modeling and computing the fracture phenomena in solids. Incremental perturbed minimization principles for a class of gradient-type dissipative materials are used to derive the perturbed governing equations. This analysis enables us to study the highly heterogeneous microstructure and monitor the uncertainty in failure mechanics. Several numerical examples are given to examine the efficiency of the proposed method.

\end{Abstract} 

\sectpa[Section1]{Introduction} 

\begin{figure}[!b]
	\centering
	{\includegraphics[width=0.85\textwidth]{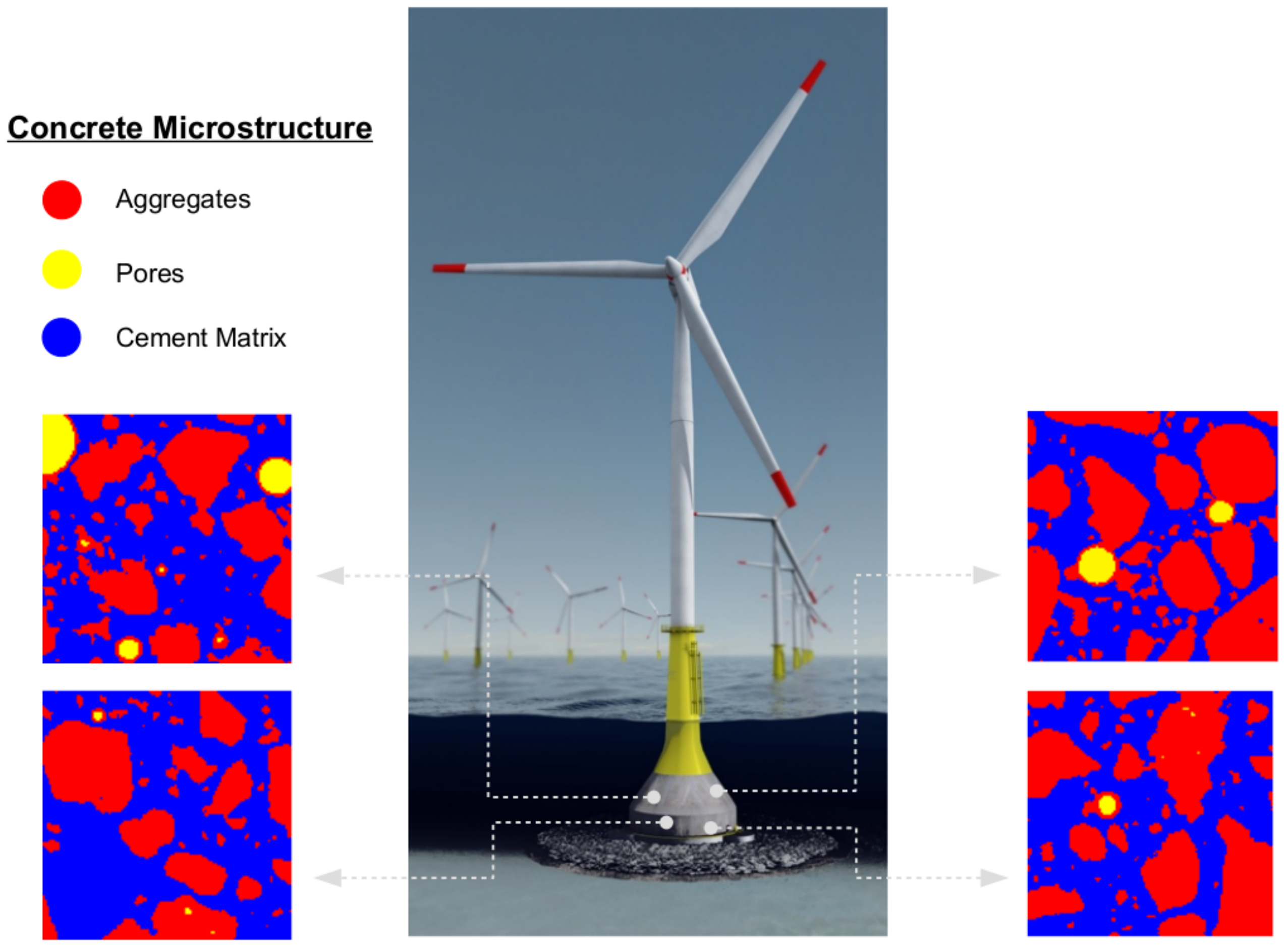}}  
	\caption{Offshore Wind Turbine (source: \texttt{germanoffshorewind.org}) with different concrete microstructures. The concrete representative volume elements (RVEs) at the microscale are consisting of aggregates, pores and cement matrix under-water (CT-images source: \texttt{www.baustoff.uni-hannover.de} related to the recent work of \cite{wessels2022computational}).}
	\label{Figure1}
\end{figure}

Investigation of crack initiation and propagation in brittle and ductile materials is a topic of intensive research to predict failure mechanisms for various engineering structures. These applications experience different failure-modes related to the desired operating conditions. Hereby, material and geometrical properties are considered to be homogeneous and deterministic at their macro-structure level. Whereas, strong fluctuation is observed in those quantities at the microstructures \cite{zohdi2004introduction}. This is quite natural as materials may contain a scatter range in their properties around a mean value. Furthermore, the well-known tolerances in the industrial manufacturing processes along with their real-life applications will produce a range of perturbations in the geometric properties. For a better understanding of the structure variation, consider the offshore wind turbine with different concrete microstructures, illustrated in Figure \ref{Figure1}. Herein, a random distribution of the aggregates and pores within the cement matrix is observed. This can {\it vary} from one point to another at the lower scale due to the {\it segmentation-tolerance} of the computer tomography CT-images, see \cite{wriggers2007micro,hain2008numerical}. 

For the safety assessment of such engineering applications, a sufficient large safety-factor is a must in the design process to account for all the uncertainties in failure mechanic problems. These applications can significantly benefit from a precisely predictive computational tool along with experimental techniques to model brittle and ductile fracture in the design phase of products. 

The computational modeling of crack propagation can be achieved in a convenient way by the continuum phase-field approach to fracture, which is based on the regularization of sharp crack discontinuities. Due to its simplicity, this methodology has gained wide interest and started to be used in the engineering community since 2008. From there on many scientists have worked in this field and developed phase-field approaches for finite elements, isogeometric analysis, and lately also for the virtual element technology. The main driving force for these developments is the possibility to handle complex fracture phenomena within numerical methods in two and three dimensions. In recent years, several {brittle \cite{BourFraMar08,KUHN10,miehe+welschinger+hofacker10,hesch+weinberg14,Wick15Adapt,rezaei2021direction,alessi2020phase,denli2020phase,heider2021review,schreiber2020phase,steinke2019phase,arash2021finite,fantoni2019phase,aldakheel2021feed,Wick+2020,pillai2020anisotropic,heider2022self,rezaei2022anisotropic,selevs2021general,ambatiphase,zhuang2022phase,seiler2021phase,tan2022phase}} and {ductile \cite{ambati+gerasimov+lorenzis15,aldakheel2020microscale,alessi2017,shanthraj2016phase,choo18,fang2019phase,kruger2019porous,nguyen2019multiscale,dean2020phase,ali2021residual,bryant2021phase,storm2021comparative,ulloa2022micromechanics,khalil2022generalised}} phase-field fracture formulations have been proposed in the literature. These studies range from {the} modeling of 2D/3D small and large strain deformations, variational formulations, multi-scale/physics problems, mathematical analysis, different decompositions, and discretization techniques with many applications in science and engineering. All these examples and the citation therein demonstrate the potential of phase-field for crack propagation.

 \begin{figure}[!b]
	\centering
	{\includegraphics[clip,trim=0.1cm 1cm 0cm 6cm, width=17.5cm]{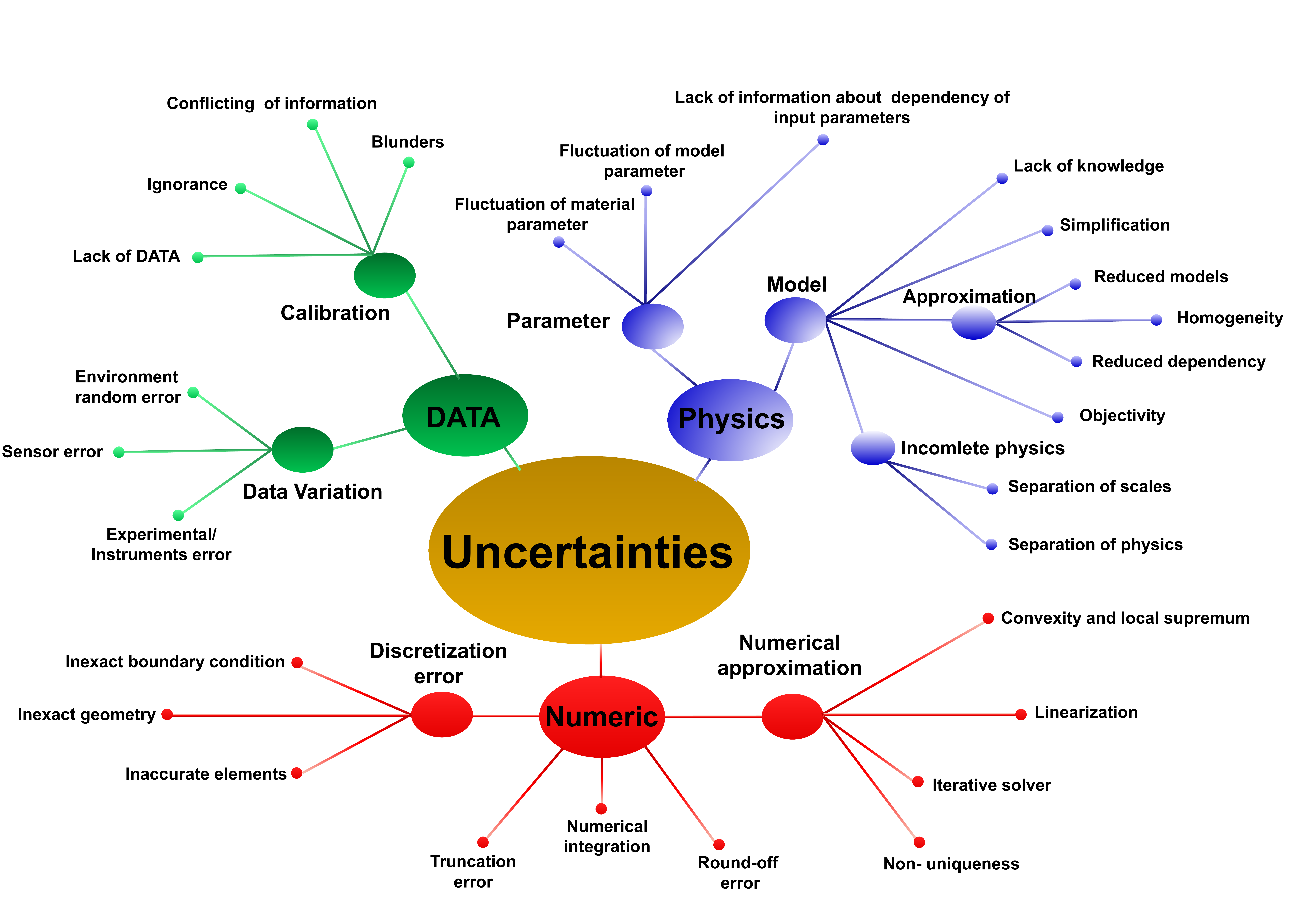}}  
	\caption{ \textcolor{black}{Sources of uncertainty in physical modeling and the numerical discrepancy.}}
	\label{Figure_source}
\end{figure}

Due to the deterministic nature of the phase-field approaches, non-unique solutions are explored for the material, geometric, and meshing perturbation. This raises the feasibility of possible several solutions and their influence on the design process. Hence, a detailed study of such randomness in those properties along with the associated system response is {\it inescapable}. To this end, we utilize the probabilistic approach to the deterministic solution, which gives us an estimation of the bounds of system response. Specifically, this work is devoted to a rigorous mathematical formulation of the stochastic-based variational framework of failure mechanisms at the micro-level. The key goal of development is to predict the failure response of materials for certain randomness and fluctuations of different phases in the highly heterogeneous microstructure. In this regard, the Monte Carlo finite element method (MC-FEM) is employed to solve the stochastic PDE-based model and approximate the expectation and the variance of the solution field of brittle/ductile failure by evaluating a large number of samples. In the MC-FEM, finite elements are utilized to discretize the computational domain and the random points according to the probability distribution to model the uncertainty \cite{guoliang1993monte,van1995use,ganesh2021quasi}. In order to improve the convergence of the random points and the computational complexity, quasi Monte Carlo techniques \cite{dick2013high,kuo2012quasi,graham2011quasi}, multilevel Monte Carlo \cite{fang2022multilevel,faustmann2021stability,kuo2017multilevel,ben2020importance} and their combination \cite{kuo2015multi,kuo2017multilevel} are proposed in literature.

Compared to the deterministic modeling, different sources of uncertainty such as fluctuations in experimental devices, noises due to the coupling of the model problems (capturing all physical phenomena from the nano-to-macro scale), and spatial variation of the material parameters lead to fluctuation in the observed results. A summary of these sources is given in Figure \ref{Figure_source}.

In the design and manufacturing process, it is essentially worthwhile to consider the material's probabilistic behavior using a microstructure.
For instance, the concrete material properties vary/fluctuate even using a similar manufacturing procedure \cite{chen2020insight,ibrahimbegovic2020reduced}. Therefore, considering the spatial variation of the concrete elastic properties, fracture energy, and plastic property (e.g., hardening) provide a more reliable modeling platform. In the Monte Carlo simulations, we generate random samples (according to the given distribution) to estimate the possible events (randomness in materials and particles spatial variations) and approximate the relative crack behavior. Of course, more number of replications will include more possible events (provides more informative data) that result in a more accurate expected value and variance. For the quasi-brittle materials, the Monte Carlo finite element method (MC-FEM) was used to model the dependence of the computed crack probabilities on the type of perturbation in \cite{gerasimov2020stochastic,ricoeur2022stochastic,colliat2007stochastic}, and the polynomial chaos expansion in functionally graded materials with random material properties is used to model the phase-field fracture, see \cite{dsouza2021non}. In computational mechanics, stochastic discretization techniques have been employed for variational theory for nonlinear problems with stochastic coefficients \cite{matthies2005galerkin,martin2022korali,pryse2017stochastic}, inelastic media under uncertainty \cite{matthies2008inelastic}, elastic-plastic material with uncertain parameters \cite{rosic2011stochastic}, fatigue crack propagation  due to the inherent uncertainties according to the material properties \cite{ben2016stochastic}, nonlinear fracture mechanics of concrete  \cite{novak2005stochastic}, and stochastic fracture response and crack growth analysis of laminated composites \cite{lal2017stochastic}. In addition to MC-FEM, different numerical methods, such as polynomial chaos expansion (PCE) \cite{beck2013stochastic}, the method of time-separated stochastic mechanics (TSM) \cite{junker2018analytical,junker2019relaxation,junker2020modeling} and stochastic finite element method  \cite{ghanem2003stochastic,stefanou2009stochastic,pryse2021neumann,reddy2008stochastic} with applications to fracture mechanics.

Recently, a Bayesian inversion approach as a probabilistic technique for the phase-field fracturing modeling has been proposed to identify material/model parameters due to the uncertainty of the fracturing material in \cite{khodadadian2020bayesian}, coupled with plasticity in \cite{noii2021bayesian}. In stochastic analysis (more specifically MC-FEM) hundreds or thousand forward runs are necessary to be performed. In numerical optimization using adjoint methods (the adjoint problem is linear, but is running backward in time) resulting in a high computational cost. Consequently, the general natural idea is to use dimension reduction techniques, as proposed in \cite{abbaszadeh2021reduced}. For reducing the computational costs of the phase-field failure analysis in a probabilistic framework (mainly Bayesian inversion), a non-intrusive global-local approach
is recently introduced, rather than using fine-scale high-fidelity finite elements \cite{noii2022bayesian}. In this type of concurrent multiscale framework, the phase-field model is solved on the fine-scale, and a linearized model (without phase-field) is employed on the global scale. While the fracture propagates, the local and global sub-domains are adjusted dynamically with the help of an adaptive predictor-corrector procedure, as shown in \cite{ALDAKHEEL2021114175}.

To explore the random nature of the material structure and its effect on the failure and fracture, this contribution first extends the prescribed model in \cite{noii2021bayesian} to a stochastic setting. The developed framework allows us to model the effect of the random distribution of the particles (densities, positions, size) and the spatial variation of the material parameters. Next, we will study the effect of the randomness on a local scale, i.e., microscopically, {in different parts of the structure} in which crack patterns can occur. The results are extended to a global approach, i.e., by computing the amount of the necessary forces for failure initiation and, therefore full fracture. Using several replications enables us to provide an accurate global pattern. Thus, by taking the expectations, the results can be extended to the whole structure. Furthermore, we will have an interval (between maximum and minimum of the forces) to determine how much force (at least) is needed for the fracture in a part of the domain and applying which amount of forces will give rise to a full fracture.

The rest of the paper is organized as follows. In Section \ref{Section2int}, a stochastic phase-field framework for modeling fracture in brittle and ductile materials will be introduced. In Section \ref{Section3stochAlloc}, an allocating strategy will be developed to model the random distribution of aggregates/pores and cement matrix in the concrete structure for two- and three-dimensional simulations. In Section \ref{Section4}, we present different multi-dimensional test experiments to model crack behavior for ductile and brittle concrete using the stochastic framework and the allocating strategy. Finally, the obtained results are summarized in Section \ref{Section5}.


\sectpa[Section2int]{Stochastic phase-field modeling of fracture}
In this section, the effect of randomness, fluctuation, and variation in phase-field fracture problems will be investigated. 

\sectpb[Section2]{Primary fields and function spaces}
We consider $\calB\subset{\mathbb{R}}^{\delta}$ be an arbitrary solid geometry, $\delta=\{2,3\}$ with a smooth  boundary $\partial\calB$. We assume Dirichlet boundary conditions on $\partial_D\calB $ and Neumann boundary conditions on $\partial_N \calB := \Gamma_N \cup \mathcal{C}$, where $\Gamma_N$  explains the outer domain boundary and {$\calC\in \mathbb{R}^{\delta-1}$} points out the fracture boundary. Furthermore, we present a probability space $(\Omega, \mathbb{A},
\mathbb{P})$, where $\Omega$ indicates the set of elementary events (the sample space), $\mathbb{A}$ is the $\sigma$-algebra of all possible
events, and $\mathbb{P}\colon \mathbb{A} \to [0,1]$ is a probability
measure. A real-valued random variable $\xi:\Omega\rightarrow\mathbb{R}$ is a set of possible events ($\Omega$), mapping the probability space to the real values. A realization $\omega=(\omega_1\ldots\omega_n)$ is given on the probability space and denote as $n$-dimensional random variable. In this work, the randomness points out the stochastic distribution of the heterogeneity (inclusions, voids), their random number (according to random distribution), the relative random radius, and fluctuation in material parameters. The randomness changes the macroscopic as well as the microscopic structure and affects the corresponding stiffness. We study this effect locally (monitoring the crack propagation pattern) and globally (the variation of the load-displacement diagram).

Denoting the event $\omega\in\Omega$, the expectation function can be defined by $\mathbb{E}[\xi]:=\displaystyle\int_{\Omega} \xi(\omega)\,\mathrm{d}\mathbb{P}$ and the variance function is defined by $\mathbb{V}[\xi]:=\mathbb{E}[\xi^2]-(\mathbb{E}[\xi])^2$. We define an inner product for set of  $(\xi,\zeta)$ as a real-valued random variable $(\xi,\zeta):\Omega\rightarrow\mathbb{R}$ for a possible events as a $\langle \xi,\zeta\rangle_{\mathcal{A}}:=\mathbb{E}\left[ \xi\,\zeta\right]$ and the following Lebesgue space of the random variables using the finite variance
\begin{align}
	\mathcal{A}:=L^2(\Omega,\mathbb{R})=\{\xi:\Omega\rightarrow\mathbb{R}:~~\|\xi\|_{\mathcal{A}}^2:=	\langle \xi,\xi\rangle_{\mathcal{A}}=\mathbb{E}[\xi^2]<\infty\}.
\end{align}
Using the above-mentioned definition, we define the the covariance operator given by $\mathcal{COV}(\xi,\zeta):=\mathbb{E}\left[\left( \xi-\mathbb{E}[\xi]\right)\left( \zeta-\mathbb{E}[\zeta]\right)  \right]= \langle \xi-\mathbb{E}[\xi],\zeta-\mathbb{E}[\zeta] \rangle_{\mathcal{A}}$. Obviously, the variance function can be defined as $\mathbb{V}(\zeta)=\mathcal{COV}(\zeta,\zeta)$. We can set the uncorrelated random variables when we have $\mathcal{COV}(\xi,\zeta)$=0.
A random variable $\xi(\boldsymbol{x};\omega):\Omega\times X\rightarrow \mathbb{R}$ is related to the spatial variable $\boldsymbol{x}\in X$ and the random variable $\xi\in\Omega$. Considering the fixed random variable $\omega\in\Omega$, $\xi(\omega,\cdot)$ indicates the deterministic cases, which is one observation in the phase-field problem (ductile/brittle).   

Considering the random variable $\omega\in \Omega$, the response of the material at point $\Bx\in\calB$ and at time $t\in \mathbb{T} = [0,T]$ can be presented by the random displacement field $\Bu(\Bx,t;\omega): \calB \times\mathbb{T}\times\Omega \rightarrow \mathbb{R}^\delta $, the random crack phase-field $d(\Bx,t;\omega): \calB \times \mathbb{T}\times\Omega \rightarrow [0,1]$.
Here, $d(x,t;\cdot)=0$ and $d(x,t;\cdot)=1$ characterize an undamaged and a completely fractured material state, respectively. The loading time interval can be discretized as
\begin{align}
	0<t_1<t_2<\cdots<t_n<\cdots <t_N=T.
\end{align}
We note that for any variable used from now onward $\bullet_n=\bullet (t_n)$. 
Having a random variable $\omega$, with the purpose of stating variational principles, we introduce the following function spaces
\begin{align}
		\boldsymbol{U}&:=\{\Bu\in\mathbf{H}^1(\calB) \quad: \quad \Bu=\overline{\Bu} \;\; \text{on} \ \partial_D\calB\},\\
		V&:=\{d\in \mathrm{H}^1(\calB{\color{black}})\quad :\quad d\geq d_n,~~d(t=0)=0\}
	\label{eq:spaces_d},
\end{align}
where $\mathbf{H}^1(\calB)=\left(H^1(\calB)\right)^{\delta}$ and $d_n$ is the damage value in a previous time instant which introduces the evolutionary character of the phase-field, incorporating an irreversibility condition in incremental form. 

In the case of von-Mises plasticity theory, we define the plastic strain tensor\\ $\Bve^p(\Bx,t;\omega): \calB \times \mathbb{T}\times\Omega \rightarrow \mathbb{R}^{\delta\times\delta}_\mathrm{dev}$ and the hardening variable $\alpha(\Bx,t;\omega): 	\calB \times \mathbb{T}\times\Omega \rightarrow \mathbb{R}_+$.
Here, $\mathbb{R}^{\delta\times\delta}_\mathrm{dev}:=\{\Be\in\mathbb{R}^{\delta\times\delta} \ \colon \ \Be^T=\Be,\ \tr{[\Be]}=0\}$ is the set of symmetric second-order tensors with vanishing trace. 
Since, gradient non-local plastic theory is employed here, the plastic strain tensor is considered as a local internal variable, while the hardening variable is a non-local internal variable. Therefore, the rate of the hardening variable $\alpha$ follows the evolution equation 
\begin{equation}
	\dot\alpha = \sqrt{\frac{2}{3}}\,\vert\dot{\Bve}^p\vert.
	\label{s2-evol-alpha}
\end{equation}
At the first time step, $\alpha(\Bx,0,\cdot) $ can be viewed as the equivalent plastic strain, which starts to evolve from the initial condition $\alpha = 0$. Concerning the function spaces, we assume sufficiently regularized plastic responses, i.e., endowed with hardening and/or non-local effects, for which we assume $\Bve^{p}\in\mathbf{Q}:=\mathrm{L}^2(\calB;\mathbb{R}^{\delta\times\delta}_\mathrm{dev})$. Moreover, in view of~\eqref{s2-evol-alpha}, it follows that $\alpha$ is irreversible. Assuming in this section the setting of gradient-extended plasticity, we define the function spaces
\begin{equation}
	Z:=\{\alpha\in H^1(\calB)\quad \colon \quad \alpha = \alpha_n + \sqrt{2/3}\,\vert z \vert, z\in \mathbf{Q},~	{\alpha (t=0)=0} \ \}.
	\label{Walpha}
\end{equation}
 The hardening law~\eqref{s2-evol-alpha} is thus enforced in incremental form by restricting the solution space $Z$ where $z\in \Bve^p-\Bve^p_n$ (shown as $Z_{\Bve^p-\Bve^p_n}$). Considering the random variable $\omega$, the gradient of the displacement field defines the symmetric strain tensor of the geometrically linear theory as
\begin{equation}
	\Bve(\Bx,t;\omega) = \frac{1}{2} [\nabla\Bu(\Bx,t;\omega) + \nabla\Bu(\Bx,t;\omega)^T]
	.
	\label{s2-disp-grad}
\end{equation}
In view of the small strain hypothesis and the isochoric nature of the plastic strains, the strain tensor is additively decomposed into an elastic part $\Bve^e$ and a plastic part $\Bve^p$ as
\begin{equation}
	\Bve(\Bx,t;\omega) = \Bve^e(\Bx,t;\omega) + \Bve^p(\Bx,t;\omega) 
	\WITH
	\tr{[\Bve]} = \tr{[\Bve^e]}
	.
	\label{s2-strain-e-p}
\end{equation}

\sectpb[Section22]{Variational principles}
\label{sec:energy_functions}

Let $\BfrakC$ denote the set of constitutive state variables. In the most general setting considered in this study, one has
\begin{equation}
	\BfrakC(\Bx,t;\omega) := \{ \Bve, \Bve^p, \alpha, d, \nabla \alpha, \nabla d \},
	\label{state}
\end{equation}
where the random variable $\omega\in\Omega$ affects the functions $(\Bu,\,\alpha,\,d)$.
In order to derive the perturbed variational formulation in the stochastic space, we set the perturbed energy density function per unit volume $ W\big(\BfrakC(\Bx;\omega);\Bfrakq(\Bx;\omega)\big)$, such that
\begin{equation}\label{bochner norm_m}
\texttt{W}:=\E\Big[W\bigg(\BfrakC(\Bx;\omega);\Bfrakq(\Bx;\omega)\bigg)\Big]=	\int_\Omega  {W}(\cdot;\omega)  \,\P(\omega),
\end{equation} 
where $\Bfrakq (\cdot)$ is a random quantity denoting the randomness in the geometry (e.g. inclusions/voids) along with the material parameters. Hereby, for a fixed random variable,  the perturbed energy function is additionally decomposed into a perturbed elastic contribution $\texttt{W}_{elas} (\cdot\,;\Bfrakq)$,  a perturbed plastic contribution $\texttt{W}_{plas} (\cdot\,;\Bfrakq)$, and a perturbed fracture contribution~$\texttt{W}_{frac} (\cdot\,;\Bfrakq)$ results in
\begin{equation}
	{W(\BfrakC;\Bfrakq):= \texttt{W}_{elas}(\Bve,\Bve^p ,d ,\alpha;\Bfrakq) + \texttt{W}_{plas}(\alpha, d, \nabla \alpha;\Bfrakq) +
		\texttt{W}_{frac}(d, \nabla d ;\Bfrakq), } 
	\label{psuedo-energy}
\end{equation}
and therefore by taking the expectation we have
\begin{align}
	{\texttt{W}=\E[W]= \E\big[\texttt{W}_{elas}(\Bve,\Bve^p ,d ,\alpha;\Bfrakq)\big] + \E\big[\texttt{W}_{plas}(\alpha, d, \nabla \alpha;\Bfrakq)\big] +\E\big[
		\texttt{W}_{frac}(d, \nabla d ;\Bfrakq)\big]. } 
	\label{psuedo-energy1}
\end{align}

%
%
\sectpb[Section221]{Computing statistical moments}
\label{sec:Computing statistical}

Let us compute statistical moments for the response observed during different noises arises from  elastic-plastic setting. To do so, we introduce the quantity of interest as $\BfrakJ(\Bu,d,\alpha;\omega)$. Considering the stochastic space, $M$ different realizations are employed to compute the expected value and afterward the variance function as
\begin{equation*}
	\begin{aligned}
		\mathbb{E}[\BfrakJ] = \E\Big[\BfrakJ(\Bu,d,\alpha;\omega)\Big]&\approx \text{E}_{\text{MC}} [\BfrakJ]:= \frac{1}{M} \sum_{i=1}^M \BfrakJ(\Bu^{(i)},d^{(i)},\alpha^{(i)})\\
	\mathbb{V}[\BfrakJ] =\mathbb{V}\Big[{\BfrakJ(\Bu,d,\alpha;\omega)}\Big]&=
	\mathbb{E}[\BfrakJ^2(\Bu,d,\alpha;\omega)]-(\mathbb{E}[\BfrakJ(\Bu,d,\alpha;\omega)])^2 
	 \\
&\approx \text{V}_{\text{MC}}[\BfrakJ]:=	 \frac{1}{M}\sum_{i=1}^M \Big(\BfrakJ(\Bu^{(i)},d^{(i)},\alpha^{(i)}) - 	 \E[\BfrakJ(\Bu^{(i)},d^{(i)},\alpha^{(i)})]\Big)^2.
		\end{aligned}
\end{equation*}
For the expected value function, we  assume two random variables $\omega$ and $\zeta$ belong to $\Omega$. Then,  $\BfrakJ(\omega)$ and $\BfrakQ(\zeta)$ are two random scalar-valued function along with  their expected scalar-valued function $\mathbb{E}[\BfrakJ]$ and $\mathbb{E}[\BfrakQ]$, respectively, so the following properties holds:
\begin{itemize}
	\item $\mathbb{E}[\BfrakJ+\BfrakQ]=\mathbb{E}[\BfrakJ]+\mathbb{E}[\BfrakQ]$ holds, if both expected values are finite.
	\item For the constant $c\in\mathbb{R}$, we have $\mathbb{E}[c\,\BfrakJ]=c\,\mathbb{E}[\BfrakJ]$.
	\item Considering the Jensen's inequality and the convex property of the norm function, we have $\|\mathbb{E}[\BfrakJ]\|\leq \mathbb{E}[\|\BfrakJ\|].$
	\item Cauchy-Schwarz inequality: $\mathbb{E}\left[\|\BfrakJ\,\BfrakQ\|\right]\leq {\displaystyle\left(\mathbb{E} \left[\|\BfrakJ\|^2]\right]\right)^{1/2}\left(\mathbb{E} \left[\|\BfrakQ\|^2\right]\right)^{1/2}}.$
	\item Minkowski inequality for the $H^1$-space:  ${\displaystyle\left(\mathbb{E}\left[\|\BfrakJ+\BfrakQ\|^2\right]\right)^{1/2}\leq \left(\mathbb{E} \left[\|\BfrakJ\|^2\right]\right)^{1/2}+\left(\mathbb{E} \left[\|\BfrakQ\|^2\right]\right)^{1/2}}.$
\end{itemize}
Regarding the spatial discretization, for a fixed time step $t\in\mathbb{T}$, and fixed random variable $\omega\in\Omega$, we assume that ${\texttt E}_h=\{E_1,\dots,E_{nel}\}$  is a quasi-uniform mesh defined in $\calB_h \approx \calB$
with mesh size 
$h := \max_{E_j\in {\texttt E}_h}{\rm diam}(E_j)$.
For the sake of simplicity, we use lowest order Galerkin discretization in 
$\mathcal{B}_{h} := \boldsymbol{S}^{1,1}_u(\texttt E_h) \times S^{1,1}_d(\texttt E_h) \times S^{1,1}_\alpha(\texttt{E}_h)$, where
\begin{align*}
	\BS ^{1,1}_u (\texttt E _h) &:= \{ \Bu \in \BH^1 (\calB) \quad : \quad \Bu | _ E \in \BP _1(E)  \quad \forall\; E \in \texttt E _ h \},\\
	S ^{1,1}_d (\texttt E _h) &:= \{ {d} \in H^1 (\calB) \quad ~: \quad {d | _ E} \in {P _1(E)}  \quad ~\,\forall\;E \in \texttt E _ h \},\\
			S ^{1,1}_\alpha (\texttt E _h) &:= \{ {\alpha} \in H^1 (\calB) \quad ~: \quad {\alpha | _ E} \in {P _1(E)}  \quad \,\,\forall\;E \in \texttt E _ h \} .
\end{align*}
Here $\BP_1(E)$ and $P_1(E)$ indicate the vectorial and scalar space of polynomials of total degrees less or equal than one, respectively \cite{hughes2012finite}. Hence, we define
\begin{align}
\BU_h&:=\left\{\Bu_h\in    \BS ^{1,1}_u (\texttt E _h)\quad : \quad \Bu_h |_{\partial_D\calB}=\overline{\Bu}  \right\},\\
V_h &:= \left\{d_h\in S ^{1,1}_d (\texttt E _h)\quad \,: \quad  d_h\geq d_{h_{n}},~~d_h(t=0)=0\right\},\\
Z_h&:=\left\{\alpha_h\in S^{1,1}_\alpha (\texttt E _h)\quad : \quad \alpha_h = \alpha_{h_{n}} + \sqrt{2/3}\,\vert z_h \vert, z_h\in \mathbf{Q},~	{\alpha_h (t=0)=0}  \right\}.
\end{align}
Then, we define the continuous solution space $\BX:=\BU\times V\times Z$ with the corresponding norm $\|\cdot\|_{\BX}$ and the discrete solution space $\BX_h:=\BU_h\times V_h\times Z_h$ which is a subset of $\BX$. Considering a fixed random variable $\omega\in \Omega$, the quantity $\BfrakJ(\Bu,d,\alpha;\omega)$ can be approximated by $\BfrakJ_h(\Bu,d,\alpha;\omega)\approx\BfrakJ(\Bu_h,d_h,\alpha_h;\omega)$.

In MC-FEM simulations, to obtain an accurate estimation of the stochastic solution, a sufficiently small mesh size in addition to several number of evaluations are needed. To this end, we define the {Bochner space $L^2(\Omega;\BX)$} for the function $\mathcal{Y} $, giving
\begin{align}
 	\| \mathcal{Y} \|_{L^2(\Omega;\BX)} :=\Big( \int_\Omega \| \mathcal{Y}(\cdot;\omega) \|_\BX^2 \,\P(\omega)
 	\Big)^{1/2}
 	= \E\Big[\| \mathcal{Y}(\cdot;\omega) \|_\BX^2 \Big]^{1/2}.
 	\label{bochner}
\end{align}
Here, the variance function is given by $\mathbb{V} ( \mathcal{Y}) = \|\E[ \mathcal{Y}]- \mathcal{Y}\|_{L^2(\Omega;\BX)}^2$. With respect to the mesh size $h$
the discretization error is computed by
\begin{align}\label{17}
\mathcal{I}_h:= {\displaystyle \left(\|\mathbb{E}\left[ \Bu-\Bu_h\right]\|_{\BX}+ \|\mathbb{E}\left[ d-d_h\right]\|_{\BX}+ \|\mathbb{E}\left[ \alpha-\alpha_h\right]\|_{\BX} \right)}.	
\end{align}
Let us assume $\chi$ is a member of ($\Bu,d,\alpha$), and $\chi_h$ is one of the approximations ($\Bu_h,d_h,\alpha_h$).  The following lemma denotes the convergence of the statistical Monte Carlo estimator. 
\begin{lemma}\label{lemma:general MC error}For the number of samples $M$, 
	$\chi \in {L^2(\Omega;\BX)}$ satisfies \cite{taghizadeh2017optimal}
	\begin{equation}\label{general MC error}
		\| {\mathbb{E}}[\chi]-\EMC[\chi] \|_{L^2(\Omega;\BX)}
		= M^{-1/2} \mathbb{V}[\chi].
	\end{equation}
\end{lemma}

\begin{proof}
	We use the defined Bochner norm in \eqref{bochner} 
	\begin{align*}
		\| {\mathbb{E}}[\chi]-\EMC[\chi] \|^2_{L^2(\Omega;\BX)}
		&=\E\Big[ \Big\| \E[\chi] - \frac{1}{M}\sum_{i=1}^{M}\chi^{(i)}  \Big\|_\BX^2\Big]\\ 
		&=\frac{1}{M^2}\sum_{i=1}^{M}\E\Big[ \| \E[\chi] - \chi^{(i)}  \|_\BX^2\Big]\\ 
		&=\frac{1}{M}\E\Big[ \| \E[\chi] - \chi  \|_\BX^2\Big]
		=M^{-1}\mathbb{V}^2[\chi].
	\end{align*}
\end{proof}
As the next step, we can compute the total error denoting the discretization error as well as the statistical error. 
Lemma \ref{lemma:general MC error} controls of the sampling
error. The discretization error also relates to used polynomial order in the  finite element method.

\begin{proposition}\label{prop2}
	For $\chi \in {L^2(\Omega;\BX)}$ and its finite element approximation $\chi_{_h}$, we assume that we have the  convergence rate $\nu$ of the discretization error \cite{khodadadian2020adaptive}
	\begin{equation}\label{MC:assump1}
	{\|\mathbb{E}[\chi-\chi_{_h}]\|_{L^2(\Omega;\BX)} \le a h^\nu,}
	\end{equation}
	and we have the upper bound for the variance estimator
	\begin{equation}\label{MC:assump2}
		\mathbb{V}[\chi_{_h}]\le b,
	\end{equation}   
where $a$ and $b$ are positive constants, introduced in \cite{khodadadian2020adaptive}.
For the MC-FEM estimator, we have the following upper error bound denoting discretization and statistical error
	\begin{equation}\label{MC error}
		\| \E[\chi]  - \EMC[\chi_h] \|_{L^2(\Omega;\BX)} \le {ah^\nu + b  M^{-1/2}}
		= O(h^\nu)+O(M^{-1/2}).
	\end{equation}
\end{proposition}

\begin{proof}
	By defining the root mean square error (RMSE), employing  the triangle inequality and Lemma~\ref{lemma:general MC
		error}, we will have 
	\begin{equation}\label{MC-error1}
		\begin{split}
			\RMSE
			&:=   \| \E[\chi]  - \EMC[\chi_h]  \|_{L^2(\Omega;\BX)} \\
			&\le  \| \E[\chi]  - \E[\chi_h]  \|_\BX +  \| \E[\chi_h]  - \EMC[\chi_h]  \|_{L^2(\Omega;\BX)}\\
			&\le \| \E[\chi-\chi_h] \|_\BX  +  M^{-1/2} \mathbb{V}[\chi_h]\\
			&\le {ah^\nu + b  M^{-1/2}}\\
			&= O(h^\nu)+O(M^{-1/2}).
		\end{split}
	\end{equation}
\end{proof}
The above-mentioned proposition points out that by reducing the mesh size and increasing number of replications, the total error reduces. Thus, following \eqref{MC:assump1},
 for the couples system of equations ($\Bu,d,\alpha$), we have 
\begin{align}\label{23}
\|\mathbb{E}\left[ \Bu-\Bu_h\right]\|_{\BX}\leq a_1 h^{\nu_{_1}}\quad  \|\mathbb{E}\left[ d-d_h\right]\|_{\BX}\leq a_2 h^{\nu_{_2}}\quad  \|\mathbb{E}\left[ \alpha-\alpha_h\right]\|_{\BX}\leq a_3 h^{\nu_{_3}}.
\end{align}
So, by replacing \eqref{23} in \eqref{17} and by defining $\hat{\nu}:=\max\{\nu_{_1},\,\nu_{_2},\,\nu_{_3}\}$, we will have
 $\mathcal{I}_h\leq \hat{a} h{^{\nu^\star}}$ where $\hat{a}$ represents the three positive constants.
%

\sectpb[Section22]{Perturbed rate-dependent functionals}
\label{sec:specifi energy}

Herein, an extension of the model proposed in~\cite{noii2021bayesian} is considered by further perturbing the energy functional and domain space (the stochastic setting). We can define a perturbed pseudo potential energy functional augmented with the random quantity {${\Bfrakq}(\boldsymbol{x};\cdot)$} as following:
\begin{equation}
	\begin{aligned}
		\calE(\Bu,\Bve^p,\alpha,d;\Bfrakq) &:= 
		 \int_{\calB} {W}(\BfrakC;\Bfrakq) \, dv  
		\; - \; 
		\vphantom{\frac{d}{dt}}
		\mathcal{E}_{ext} (\Bu;\Bfrakq), 
		\label{potential-functional}
	\end{aligned}
\end{equation}
here $\mathcal{E}_{ext}$ considers the perturbed external loads as
\begin{equation}
	\mathcal{E}_{ext} (\Bu;\Bfrakq) := 
	  \int_{\calB} \overline\Bf \cdot \Bu\, dv\,  +
 \int_{\partial_N\calB} \overline{\Btau} \cdot \Bu\, da.
\end{equation}
The quantity $\Bfrakq$ indicates the perturbered material parameter in the homogeneous and heterogeneous structures can be defined using the perturbation value $\eta$ as
\begin{align}
	\Bfrakq=\bar{\Bfrakq}+\eta \tilde{\Theta},
	\label{theta}
\end{align}
where $\tilde{\Theta}$ is a uniformly distributed random variable in [-1,\,1], and $\eta$ denotes as the material parameters variations. For the phase-field fracture in brittle and ductile materials, a set of parameters, i.e., $\bar{\Bfrakq}\in\{E,\mu,K,G_c,\psi_c,H,\sigma_Y\}$ is given in Table \ref{table1}. Denoting the variation parameter $\eta$, the perturbed materials can be defined point-wise (heterogenous), or for the whole domain (homogeneous). We use this notation to point out the fluctuation in the material property (homogeneous/heterogeneous) in addition to the randomness due to the random distribution of the particles (aggregates/voids). 

Considering the effect of the randomness, the energy functional is defined as
\begin{align}
\calE(\Bu,\Bve^p,\alpha,d)=
 \E\Big[\int_{\calB} {W(\BfrakC;\Bfrakq)} \, dv\Big] 
\; - \; 
\vphantom{\frac{d}{dt}}
\E\left[\mathcal{E}_{ext}  (\Bu;\Bfrakq) \right].
\end{align}
%
Next, to formulate the variational formulation setting, it is required to define the perturbed constitutive energy density functions, namely $\texttt{W}_{elas}$, $\texttt{W}_{plas}$, and $\texttt{W}_{frac}$. 

\sectpc{Elastic energy contribution}
The elastic energy density function $\texttt{W}_{elas}$ in \eqref{psuedo-energy} formulated based on the effective strain energy density $\psi_e(\Bve^e;\Bfrakq)$. Here, the perturbed effective strain energy density function is additively decomposed into \textit{fracturing} and \textit{unfracturing} parts is employed. Thus, the strain tensor is decomposed into \textit{volume-changing}
(volumetric) and \textit{volume-preserving} (deviatoric) to avoid failure in compression parts, as
\[
\bm\varepsilon^e(\Bu,\Bx;\omega)=\bm\varepsilon^{e,vol}(\Bu,\Bx;\omega)+\bm\varepsilon^{e,dev}(\Bu,\Bx;\omega),
\]
where
\begin{equation}
	\begin{aligned}
	&\bm\varepsilon^{e,vol}(\Bu,\Bx;\omega):=\frac{1}{3}(\bm \varepsilon^e((\Bu,\Bx;\omega):\text{\BI}){\text{\BI}},\\
	&\bm\varepsilon^{e,dev}(\Bu,\Bx;\omega):=\mathbb{P}:\bm \varepsilon^e,
	\WITH \mathbb{P}:=\mathbb{I}-\frac{1}{3}\text{\BI}\otimes\text{\BI} \AND 
	\mathbb{I}_{ijkl}:=\frac{1}{2}\big(\delta_{ik}\delta_{jl}+\delta_{il}\delta_{jk}\big).
	\end{aligned}
\end{equation}
The perturbed effective strain energy function $\psi_e(\Bve^e;\Bfrakq)$ reads:
\begin{equation}
	\psi_e\Big(\text{I}_1(\Bve^e;\omega),\text{I}_2(\Bve^e;\omega);\Bfrakq(\Bx;\omega)\Big)=\psi_e^{+}\Big(\text{I}_1,\text{I}_2;\Bfrakq\Big)+\psi_e^{-}\Big(\text{I}_1,\text{I}_2;\Bfrakq \Big),
\end{equation}
such that
\begin{equation}
	{\psi_e^{+}}={\text{H}{^+}[\text{I}_1]}\psi_e^{vol}\big(\text{I}_1;\Bfrakq\big)
	+\psi_e^{dev}\big(\text{I}_1,\text{I}_2;\Bfrakq\big)~\AND 
	{\psi_e^{-}}=\big(1-{\text{H}{^+}[\text{I}_1]}\big)\psi_e^{vol}\big(\text{I}_1;\Bfrakq\big).~	
	\label{eq232_5}
\end{equation}
{Therein, $\text{H}{^+}[\text{I}_1(\Bve^e;\omega)]$ is a \textit{positive Heaviside function} which returns one and zero for $\text{I}_1(\Bve^e;\omega)>0$ and $\text{I}_1(\Bve^e;\omega)\leq0$, respectively.} We note that the volumetric and deviatoric counterpart of energy admits following additive split:
\begin{equation}
	\psi_e^{vol}\big(\text{I}_1;\Bfrakq\big)
	=\frac{K}{2}\text{I}^2_1
     \AND
	\psi_e^{dev}\big(\text{I}_1,\text{I}_2;\Bfrakq\big)
	=\mu\Big(\frac{\text{I}_1^2}{3}-\text{I}_2\Big),
	\label{eq:psi_iso11}
\end{equation}
in terms of the the bulk $K$ and shear modulus $\mu$, where $\text{I}_1(\Bve^e;\omega):=\text{tr}[\Bve^e]$  and $\text{I}_2(\Bve^e;\omega):=\tr[(\Bve^e)^2]$ denote the first and second invariants. So the total elastic contribution to the pseudo-energy \eqref{psuedo-energy} finally reads
\begin{equation}
	\texttt{W}_{elas}(\Bve,\Bve^p,d,\alpha;\Bfrakq):=
	g(d)\; \psi_e^{+}(\text{I}_1,\text{I}_2;\Bfrakq)+ \psi_e^{-}(\text{I}_1,\text{I}_2;\Bfrakq),
	\label{elas-part0}
\end{equation}
such that
\begin{equation*}
 \E\Big[	\texttt{W}_{elas}(\Bve,\Bve^p,d,\alpha;\Bfrakq)\Big]\approx \text{E}_{\text{MC}}\Big[ 	\texttt{W}_{elas}(\Bve,\Bve^p,d,\alpha;\Bfrakq) \Big],
	\label{elas-part}
	\tag{$ {\text{E}}$}
\end{equation*}
where $+g(d(\Bx;\Bfrakq))$ is the \textit{degradation function}.

\sectpc{Fracture energy contribution} 
The phase-field contribution $\texttt{W}_{frac}$ is expressed in terms of the crack surface energy density $\gamma_l$ and the regularized fracture length-scale parameter $l_f$  to smooth fracture sharp response. In favor of regularization, following ~\cite{miehe+welschinger+hofacker10}, the sharp-crack surface topology $\calC$ is modified by a smooth functional $\calC_l$. The regularized functional reads
\begin{equation}
	\calC_l(d) = \E\Big[\int_{\calB} \gamma_l(d, \nabla d) \, dv\Big] .
	\label{s2-gamma_l}
\end{equation}
For $\omega\in \Omega$, the standard density function for the $\gamma_l$ is defined as 
\begin{equation}
	\gamma_l(d, \nabla d;\Bfrakq):=\frac{1}{c_f}\, \bigg(\frac{\Delta(d)}{l_f} + l_f \nabla d \cdot \nabla d \bigg) \WITH {c_f:=4\int_0^1\sqrt{\Delta(b)}\,db ,}
\end{equation}
where $\Delta(d)$ is a monotonic and continuous \emph{local fracture energy function} such that $\Delta(0)=0$ and  $\Delta(1)=1$ where the effect of the randomness is taken into account as well. In the following, two different accepted formulation denoted as a linear (with elastic stage) and quadratic (without elastic stage) order are formulated. Hence, for a fixed event $\omega$ we define 
\begin{equation}
	\Delta(d,\cdot) : = \begin{cases}  d\phantom{^2} \implies c_f=8/3 \quad &\texttt{AT-1}, \\ d^2 \implies c_f=2 \quad &\texttt{AT-2}.  \end{cases}
	\label{wd}
\end{equation}
Thus, the perturbed fracture contribution related to ~\eqref{psuedo-energy} is computed denoting the randomness, i.e., 
\begin{equation}
	\begin{aligned}
		\texttt{W}_{frac}(d, \nabla d;\Bfrakq):=
		g_f \gamma_l(d, \nabla d;\Bfrakq)\Big]  ,
	\end{aligned}
	\label{frac-part0}
\end{equation}
such that
\begin{equation*}
\E\Big[\texttt{W}_{frac}(d, \nabla d;\Bfrakq)\Big] 	\approx \text{E}_{\text{MC}}\Big[ 	\texttt{W}_{frac}(d, \nabla d;\Bfrakq)\Big] ,
	\label{frac-part}
	\tag{$ {\text{F}}$}
\end{equation*}
where $g_f$ is a parameter that allows to recover different models. This will be formulated in Section 2.5.

\sectpc{Plastic energy contribution}
The plastic energy counterpart $\texttt{W}_{plas}$ is formulated based on an effective plastic energy density denoted as $\psi_p$ in $\omega\in\Omega$ for gradient-extended von Mises plasticity, as:   
\begin{equation}
	{\psi}_{p}(\alpha,\nabla\alpha;\Bfrakq) :=\frac{1}{2}H \alpha^2(\Bx;\omega) + \frac{1}{2}{\sigma_Y}\,l_p^2\nabla\alpha(\Bx;\omega)\cdot\nabla\alpha(\Bx;\omega) ,
	\label{psi_p}
\end{equation}
here, $\sigma_Y$ is the initial yield stress, $H\ge 0$ is the isotropic hardening modulus and $l_p$ is the \emph{plastic length-scale}.
Thus, the perturbed plastic pseudo-energy density~\eqref{psuedo-energy} formulated as:
\begin{equation}
	\begin{aligned}
		\texttt{W}_{plas}(\alpha,d ,\nabla \alpha;\Bfrakq):=
		g(d) \; {\psi}_{p}(\alpha,\nabla\alpha;\Bfrakq),
	\end{aligned}
	\label{plas-part0}
\end{equation}
such that
\begin{equation*}
	  \E\Big[\texttt{W}_{plas}(\alpha,d ,\nabla \alpha;\Bfrakq)\Big]	\approx \text{E}_{\text{MC}}\Big[ 		\texttt{W}_{plas}(\alpha,d ,\nabla \alpha;\Bfrakq)\Big].
\tag{$ {\text{P}}$}
	\label{plas-part}
\end{equation*}
%

\sectpc[Section232]{Plastic dissipation}

Next, we define the plastic dissipation-potential density function. This thermodynamically consistent function provides a major restriction on constitutive equations for elastic-plastic and dissipative materials based on the principle of maximum dissipation. This thermodynamical restriction is due to the second law of thermodynamics (Clausius-
Planck inequality) within \textit{a reversible (elastic) domain} in the space of the dissipative forces.
So, let us define dual driving force $\{\Bs^p , -h^p\}$ with respect to the primary fields $\{\Bve^p,\alpha\}$. 
Following the classical von-Mises plasticity setting, the yield surface function reads
\begin{equation}
	\beta^p(\Bs^p,h^p,d;\Bfrakq):= \hbox{$\sqrt{3/2}$}\;\vert \BF^p \vert - h^p- g_p(d)\sigma_Y 
	\WITH
	\BF^p:= \dev[\Bs^p] = \Bs^p - \frac{1}{3} \mbox{tr} [\Bs^p] \text{\BI} .
\end{equation}
Thus, with the yield function at hand, dissipation-potential density function for plastic response reads
\begin{equation}
	\widehat{{\Phi}}_{plast}(\dot\Bve^p,\dot\alpha,d;\Bfrakq) = \sup_{\{\Bs^p,h^p\}} \{ \Bs^p:\dot{\Bve}^p - h^p\dot\alpha  \mid \ \beta^p(\Bs^p,h^p;d;\Bfrakq)\leq 0 \} ,
	\label{phiP}
\end{equation}
which follows from the principle of maximum plastic dissipation. Taking supremum of inequality  function \req{phiP}  yields, as a necessary condition, the primal representation of the plasticity evolution problem in the form of a Biot-type equation:
\begin{equation}
	\{\Bs^p , -h^p\}\in\partial_{\{\dot\Bve^p,\dot\alpha\}}\,\widehat{{\Phi}}_{plast}(\dot\Bve^p,\dot\alpha;d;\Bfrakq) .
	\label{eq:pbiot}
\end{equation}
Considering the effect of the randomness, the dissipation potential functional for plastic flow defined as
\begin{equation}
		\begin{aligned}
	{\mathcal{D}}_{plast}(\dot\Bve^p,\dot\alpha,d;\Bfrakq)
	=\int_\Omega \int_\calB{\widehat{{\Phi}}_{plast}}~\mathrm{d}v\,\P(\omega) =\E \Big[ \int_\calB{	\widehat{{\Phi}}_{plast}}~\mathrm{d}v\Big].
 	\end{aligned}
\end{equation}
Following \citep{miehe+welschinger+hofacker10}, the dissipative function $\widehat{{\Phi}}_{vis}$ due to viscous resistance forces is defined as
\begin{equation}
	\widehat{{\Phi}}_{vis} (\dot{d},\dot{\alpha};\omega) := \frac{\eta_f}{2}\dot{d}^{\,2} (\Bx;\omega)+I_+(\dot d)+ \frac{\eta_p}{2}\dot{\alpha}^2(\Bx;\omega)+I_+(\dot \alpha).
	\label{Dvis}
\end{equation}
Here, $\eta_f$ and $\eta_p$ are material parameters that characterize the viscous response of the fracture and plasticity evolution, respectively. In this work, we assumed that the values for $\eta_f$ and $\eta_p$ are not perturbed, so deterministic values are used. As before, a global rate potential of the dissipative power density with viscous regularized evolution reads
\begin{equation}
	\mathcal{D}_{vis}(\dot{d},\dot{\alpha}) := 
	\E\Big[\int_\calB \widehat{{\Phi}}_{vis} (\dot{d},\dot{\alpha};\Bfrakq) \, dv\Big].   
	\label{potential-functional-visco}
\end{equation}
The total dissipation potential given by:
\begin{equation*}
	\mathcal{D}(\dot{\BfrakC};\Bfrakq)=\mathcal{D}_{plast}(\dot{\BfrakC};\Bfrakq)+\mathcal{D}_{vis}(\dot{\BfrakC};\Bfrakq).
	\tag{$ {\text{Diss}}$}
\end{equation*}

\sectpc[Section232]{Minimization principle for the perturbed evolution problem}
%
Here, the governing equations of the failure analysis for brittle and ductile materials can be derived from basis of the energy functional~\eqref{potential-functional} by invoking rate-type variational principles~\citep{miehe2011}. Hence, the energy functional for the fracturing elastic-plastic solid material is required for the following potential
\begin{align}
	\Pi ( \dot{\Bve},\dot{\Bve^p}, d,\dot\alpha;\Bfrakq)
	:&=\frac{d}{dt}\widehat{\mathcal{E}} ( {\Bve},{\Bve^p},d,\alpha;\Bfrakq)
	+\calD( \dot{\Bve},\dot{\Bve^p}, d,\dot\alpha;\Bfrakq)
	-\mathcal{E}_{ext}({\dot{\bm u}};\Bfrakq),
	\label{eq:energy_rate}
\end{align}
to be minimized. Thus perturbed rate-dependent gradient-extended energy functional is miminzied through following compact form
\begin{equation*}
	\boxed{ (\dot{\Bve},\dot{\Bve^p}, d,\dot\alpha;\Bfrakq) = \arg\big\{ 
		\inf_{\Bu\in\BU} \ 
		\inf_{ d\in V} \	
		\inf_{\{\Bve^p,\alpha\}\in\mathbf{Q}\times Z_{\Bve^p-\Bve^p_n}} \,  \Pi(\dot{\Bve},\dot{\Bve^p}, d,\dot\alpha;\Bfrakq) \big\} , }
	\label{minimzation_rate_energy1}
	\tag{$ {\text{min.rate}}$}
\end{equation*}
%
\sectpb[Section23]{Perturbed incremental functional}
\label{sec:Inc_specifi energy}

In this section, to formulate of \textit{transition rules} from intact region to the fully damaged bulk response, degradation function is introduced. Specifically, the fracture phase-field enters as a geometric internal variable for both elastic and plastic contribution in a simple quadratic form through following degradation function:
\begin{equation*}
  g(d)=\big(1-d(\Bfrakq(\Bx;\omega))\big)^2
\end{equation*}
along with fracture constant $g_f=2l_f c_f \psi_c$, where $\psi_c$ is a specific critical fracture energy.

Next, to further extend a global rate potential form $\Pi(\Bu,\Bve^p,\alpha,d;\Bfrakq)$ given in \req{minimzation_rate_energy1}, in-line with our recent study in \cite{noii2021bayesian}, a perturbed incremental energy minimization based on $\Pi^*(\Bu,\Bve^p,\alpha,d;\Bfrakq)$ is defined on the finite time increment $[t_{n}, t_{n+1}]$, through following potential		
\begin{equation}\label{final_inc_pot0}
	\begin{aligned}
	\Pi^*( {\Bve},{\Bve^p}, d,\alpha;\Bfrakq)
	:&=\int^{t_{n+1}}_{t_{n}} \Pi ( \dot{\Bve},\dot{\Bve^p}, d,\dot\alpha;\Bfrakq) ~\mathrm{d}t\\
	&=\widehat{\mathcal{E}}^*( {\Bve},{\Bve^p},d,\alpha;\Bfrakq)
	+\calD^*( \dot{\Bve},\dot{\Bve^p}, d,\dot\alpha;\Bfrakq)
	+\mathcal{E}^*_{ext}(\bm u;\Bfrakq).\\
\end{aligned}
\end{equation}	
Considering the randomness, we take the expectation from the energy function as  
	\begin{equation}\label{final_inc_pot}
		\begin{aligned}
	\Pi^*( {\Bve},{\Bve^p}, d,\alpha;\Bfrakq)	&=\E\Big[\int_{\calB} \Big({W}(\BfrakC;\Bfrakq)-{W}(\BfrakC_n;\Bfrakq)\Big)~\mathrm{d}v\Big]\\
	&+\E\Big[\int_\calB
	{\tau\widehat{{\Phi}}^*_{vis}}
	+ I_+( d-d_n)\Big)~\mathrm{d}v\Big]\\
	&+\E\Big[\int_\calB{	\Big(\tau\widehat{{\Phi}}^*_{plast}}+ I_+( \alpha-\alpha_n)\Big)~\mathrm{d}v\Big]\\
	&+\E\Big[\int_{\calB} \overline\Bf \cdot ({\Bu-\Bu_n})~\mathrm{d}v
	+\int_{\partial_N\calB} \overline{\Btau} \cdot ({\Bu-\Bu_n})~\mathrm{d}a\Big],\\
\end{aligned}
\end{equation}
with
\begin{equation}
  \widehat{{\Phi}}^*_{vis}= \frac{\eta_f}{2\tau^2}(d-d_n)^{\,2} +I_+(d-d_n)+ \frac{\eta_p}{2\tau^2}({\alpha-\alpha_n}^2)+I_+(\alpha-\alpha_n),
	\label{phiP_tau123}
\end{equation}
and incremental plastic dissipation potential as
\begin{equation}
	\widehat{{\Phi}}^*_{plast}(\Bve^p,\alpha,d;\Bve_n^p,\alpha_n;\Bfrakq) = 
	\frac{1}{\tau}\sup_{\{\Bs^p,h^p\}} \{ \Bs^p:\big({\Bve}^p-{\Bve}_n^p\big) - h^p\big(\alpha-\alpha_n\big)  \ \mid \ \beta^p(\Bs^p,h^p,d)\leq 0 \} ,
	\label{phiP_tau124}
\end{equation}
where $\tau = t_{n+1} - t_n > 0$ denotes the step length. Here, both plasticity and the phase-field functions are a function of space $\Bx\in \calB$ and $\omega\in\Omega$. 
Through the incremental potential given in \req{final_inc_pot} at hand, the time-discrete counterpart of the canonical rate-dependent variational principle in \req{minimzation_rate_energy1} takes the following compact form:
\begin{equation*}
	\boxed{ ({\Bve},{\Bve^p}, d,\alpha) = \arg\big\{ 
		\inf_{\Bu\in\BU}  \	
		\inf_{d\in V} \
		\inf_{\{\Bve^p,\alpha\}\in\mathbf{Q}\times Z_{\Bve^p-\Bve^p_n}}  \, \Pi^*({\Bve},{\Bve^p}, d,\alpha;\Bfrakq)  \big\} . }
	\label{minimiaz2}
		\tag{$ {\text{min.incr}}$}
\end{equation*}
Accordingly, the global primary fields are determined through the stationarity conditions of the minimization problem~\eqref{minimiaz2}:  find $\Bu\in\BU$, $\alpha\in Z$, and $d\in V$ such that
\begin{align*}
	\left\{
	\begin{aligned}
		&\mathbb{E} \Big[
		\int_\calB\big[\Bsigma(\Bve,\Bve^p,d;\Bfrakq):\Bve(\delta\Bu)-\overline{\Bf}\cdot\delta\Bu\big]\,dv\Big]\, -\mathbb{E} \Big[
		\int_{\partial_N\calB}\overline{\Btau}\cdot\delta\Bu\,da\Big] =0 \quad  &&\forall \,\delta\Bu\in\BU,  \\
			&\mathbb{E}\Big[ 	\int_\calB \bigg(-\sqrt\frac{3}{2}\vert\BF^p(\Bve,\Bve^p,d)\vert  +  (1-d)^2\sigma_Y  + \partial_\alpha I_+(\alpha-\alpha_n)\\  
			 & \hspace*{2cm}+(1-d^2)\,H\alpha + \sigma_Y \, l_p^2(1-d)^2\nabla\alpha\cdot\nabla(\delta\alpha)\bigg)\delta\alpha\,dv\Big]  \ni 0  \quad &&\forall \, \delta \alpha\in Z. \\
		&\mathbb{E}\Big[
		\int_\calB \bigg( \big(\,(1-d)\calH -  {d}\, \big)\delta d - \frac{\eta_f}{\tau} (d-d_n)\delta d - l_f^2  \nabla d \cdot \nabla (\delta d) \bigg)\,dv\Big]= 0   \quad  &&\forall \delta d\in V,
	\end{aligned}
	\right.
	\label{phf_weak_m3_hist}
	\tag{$ {\text{M}}$}
\end{align*}
In \req{phf_weak_m3_hist} the crack driving force function shown as  $\calH(\Bx,t;\Bfrakq)$ is used to impose the damage irreversibility condition through history field as:
\begin{equation*}
	\calH(\Bx,t;\Bfrakq):=\max_{s\in[0,t]}\widetilde{D}\big(\BfrakC(\Bx,s);\Bfrakq\big) \WITH  \widetilde{D}:=\zeta{\Big\langle}\frac{\psi_e^{+} + \psi_p}{\psi_c} - 1 {\Big\rangle}\;,
	\label{histfield_m3}
		\tag{$ {\text{H}}$}
\end{equation*}
where, the Macaulay bracket denotes the ramp function $\langle x \rangle:=(x+|x|)/2$. Additionally, $\zeta\geq0$ is a scaling parameter to further providing relaxation of the formulation, allowing to tune the post-critical range \cite{noii2021bayesian}.
\begin{Remark}
	\label{rem:redu_brit_frac}
	So far, we studied ductile phase-field fracture in a stochastic space which is either elastic-plastic response followed by damage (hereafter $\texttt{E-P-D}$); or elastic, followed by elastic-plastic, and then plastic-damage (hereafter $\texttt{E-P-DP}$). To reduce the model into a brittle fracture response, it is sufficient that $\sigma_Y\rightarrow\infty$ and $l_p\rightarrow0$ (hereafter $\texttt{E-D}$) to be imposed. Additionally, in the case of $\texttt{E-D}$, we have used $\texttt{AT-2}$ in \req{wd}, while for $\texttt{E-P-D}$ $\texttt{AT-1}$ is used.
\end{Remark}

\begin{figure}[!b]
	\centering
		\includegraphics[width=\textwidth]{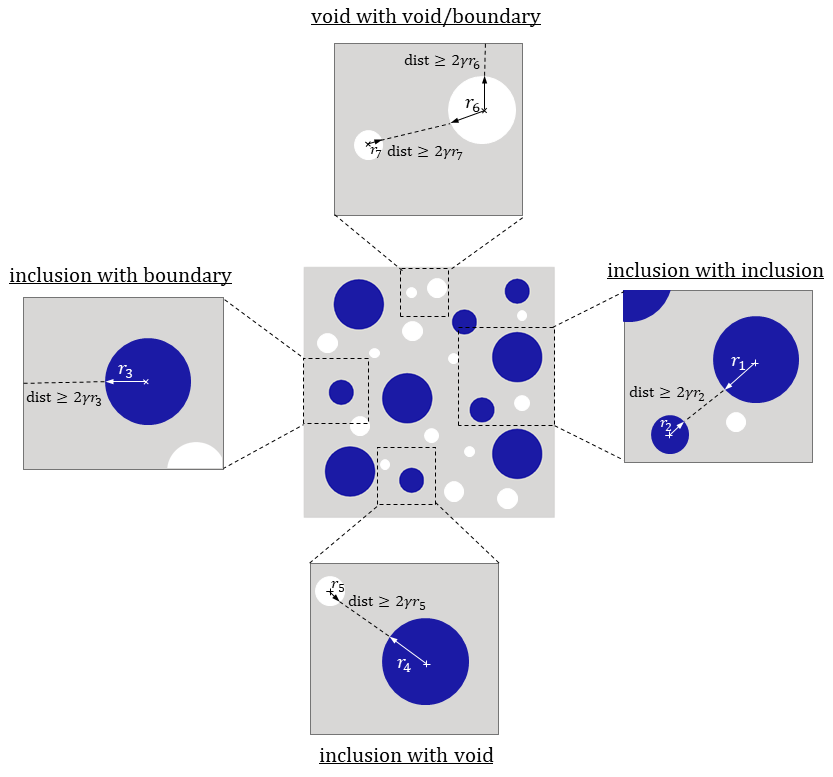}
	\caption{Allocating process of inclusions (aggregates) and voids (pores) to avoid overlapping.}
	\label{placing}
\end{figure}
%

\sectpa[Section3stochAlloc]{Random allocation of the heterogeneities}


In this section, we developed a placing strategy to define the random position, number, and size (radius) of the heterogeneities (particles) within a matrix. A uniform probability distribution is assumed for those heterogeneities (namely voids and inclusions). To this end, a given particle density is firstly considered, then we follow the procedure described in Algorithm  \ref{Algorithm1} as long as the volumes/areas of both components(voids "$\texttt V$" and inclusions "$\texttt I$")  are less than the given density. Furthermore, we should guarantee that there is no overlap between the particles and that all of them are fully allocated in the given domain (square/cube). Therefore, for each particle (void or inclusion), we consider the minimum and maximum coordinates (${\Bx}_{\text{min}}/{\Bx}_{\text{max}}$) and the radius (${\Large r}_{\text{min}}/{\Large r}_{\text{max}}$) in the Cartesian coordinate system as
\begin{equation}
	\hspace{-0.2cm}	{\Bx}={\Bx}_{\text{min}}+\theta\left( {\Bx}_{\text{max}}-{\Bx}_{\text{min}}\right) \AND
	\hspace{-0.2cm}		{\large r}={\Large r}_{\text{min}}+\theta\left( {\Huge r}_{\text{max}}-{\Large r}_{\text{min}}\right),
\end{equation}
where the random variable $\theta$ is uniformly distributed between $0$ and $1$. 
In such a heterogeneous structure, the thickness of the matrix-material can be related to the particles sizes as well documented in \cite{wittmann1985simulation,wang1999mesoscopic,wriggers2006mesoscale}. Hereby, the allowed minimum thickness of the matrix-material: a) between the heterogeneities (with radius $r_i$) and the boundary is considered as $2 \gamma \cdot r_i$ and b) between two heterogeneities is assumed as $\gamma$ times the size of the component, as shown in Figure \ref{placing}. The value of the distribution parameter $\gamma$ depends on the volume fraction of heterogeneity. To ensure this condition, we enlarge the particles size to $(1+\gamma)2r$, then we follow the placing strategy to estimate the positions; however, the radii will be determined without the enlargement. In the numerical examples, we set $\gamma=0.1$. However, the allocation strategy can be used for different values of $\gamma$, even $\gamma=0$. In summary, the chosen algorithm for the random distribution of the heterogeneities is given in Algorithm  \ref{Algorithm1}. 
\begin{figure} [!b]
	\centering
	{\includegraphics[width=1\textwidth]{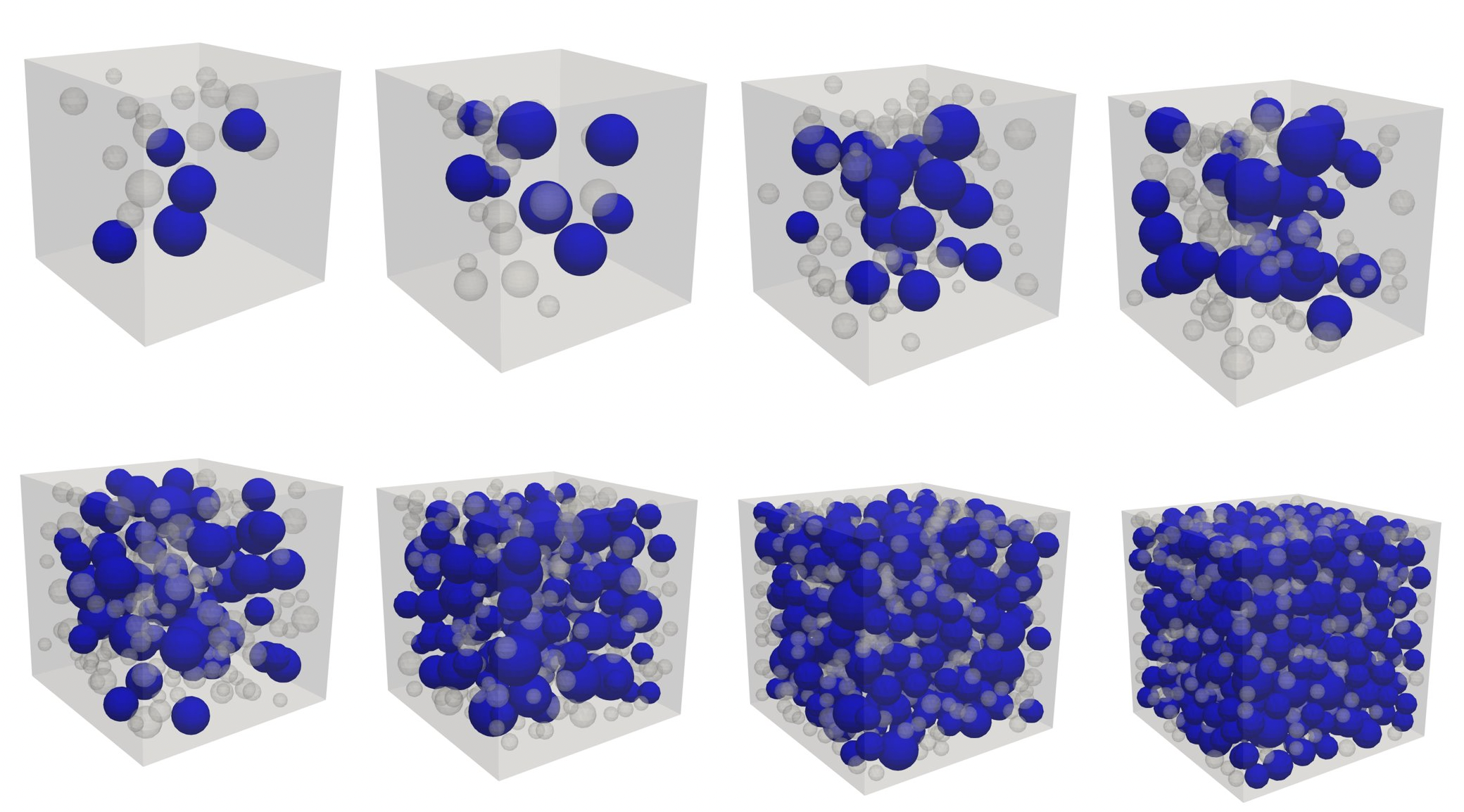}}  
	\caption{The RVE structure with different percentage of inclusions and voids.  In the first row, from left to right, the inclusions and voids densities are (3\%+3\%),  (5\%+5\%), (10\%+6\%), and (15\%+7\%). In the second row, the inclusions and voids densities are (20\%+8\%),  (30\%+10\%), (40\%+10\%), and (50\%+12\%). }
	\label{Figure1T}
\end{figure}
\begin{algorithm}
	\textbf{Inputs:}  
	$(\Bx^{\bullet}_{min},\Bx^{\bullet}_{max},r^{\bullet}_{min},r^{\bullet}_{max})$\;\text{with} $\bullet\in\{\texttt V,\texttt I \}$, voids/inclusions percentage: $\widehat{\texttt V}$/$\widehat{\texttt I}$\\[2mm] 
	\hspace*{1.6cm} 
	domain $\calB$ volume (\texttt{volume}), enlargement factor: $\gamma$.\\[4mm]
	\textbf{Initialization:}\quad$\texttt{V}=0,$\quad$\texttt{I}=0,$\\[3mm]
	\hspace*{4cm}$\bar{\texttt V}=\texttt{volume}\times\widehat{\texttt V}$\AND$\bar{\texttt{I}}=\texttt{volume}\times\widehat{\texttt{I}}$\\[4mm]
	\textbf{while} ${\texttt V\le\bar{\texttt V}}$~~\AND $\texttt{I}\le\bar{\texttt{I}}$  			
	{	\begin{equation*}
			\quad	\left\{
			\begin{array}{ll}
				\hspace{-0.2cm}	{\Bx}{^{^\texttt V}}={\Bx}_{\text{min}}^{^{\texttt V}}+\theta\left( {\Bx}_{\text{max}}^{^{\texttt V}}-{\Bx}_{\text{min}}^{^{\texttt V}}\right) \\[4mm]
				\hspace{-0.2cm}		{\large r}{^{^\texttt V}}={\Large r}_{\text{min}}^{^{\texttt{V}}}+\theta\left( {\Huge r}_{\text{max}}^{^{\texttt V}}-{\Large r}_{\text{min}}^{^{\texttt V}}\right)
			\end{array},
			\right.
			~~~	\left\{
			\begin{array}{ll}
				\hspace{-0.2cm}{\Bx}{^{^\texttt{I}}}={ \Bx}_{\text{max}}^{^{\texttt{I}}}+\theta\left( {\Bx}_{\text{max}}^{^{\texttt{I}}}-{\Bx}_{\text{min}}^{^{\texttt{I}}}\right) \\[4mm]
				\hspace{-0.2cm}{\Large r}{^{^\texttt{I}}}={\Large r}_{\text{min}}^{^{\texttt{I}}}+\theta\left( {\Large r}_{\text{max}}^{^{\texttt{I}}}-{\Large r}_{\text{min}}^{^{\texttt{I}}}\right)
			\end{array}.
			\right.
	\end{equation*}}
	
	\hspace{1cm} 1) Assume the enlarged aggregates as $2(1+\gamma)	{\Large r}{^{^\texttt I}}$.\\
	
	\hspace{1cm}	2) No overlap between the voids and aggregates.\\
	
	\hspace{1cm}	 3) The particles are within the domain boundary.\\
	
	\hspace{1cm}	 4) Determine \textit{total} fraction of voids and inclusions: \\

	$ \hspace{3cm}\qquad\texttt V=\texttt V+\frac{4}{3}\pi({\large r}{^{^\texttt V}})^3\AND\texttt{I}=\texttt{I}+\frac{4}{3}\pi({\large r}{^{^\texttt{I}}})^3$.\\
	
	\hspace{1cm} 5) Checking step
	$\texttt I+\texttt V<\displaystyle\int_{\calB}d\Bx$=\texttt{ volume} \\
	\textbf{end while}\\[4mm]
	\textbf{Outputs:} Cartesian coordinate: $\Bx=[\Bx^\texttt V,\Bx^\texttt I]^T$, radii: $r=[r^\texttt V,r^\texttt I]^T$.   
	\caption{The allocating strategy} \label{Algorithm1}
\end{algorithm}
We use the allocating strategy for the two- and three-dimensional cases based on the explained randomness and the given densities. To study the efficiency and the accuracy of the allocating algorithm, we produce RVE structures for different densities of inclusions and voids as shown in Figure \ref{Figure1T}. For the 2D distributions, we will use quadrilateral meshes and for the 3D ones, the tetrahedral meshes will be employed.

 \begin{table}
 	\caption{The material parameters for the one-, two- and three-dimensional cases. }
 	\vspace{1mm}
 	\centering
 	\begin{tabular}{|l |c c c c c cc|}
 		\hline
 		Example \hspace{0.15cm}&$E\,[\text{MPa}] $        & $\mu\,[\text{MPa}]$    &   $K\,[\text{MPa}]$   &   $G_c\,[\text{MPa.m}]$   & $\psi_c\,[\text{MPa}]$  & $H\,[\text{MPa}]$ & $\sigma_Y\,[\text{MPa}]$  \\\hline
 		1D-brittle \hspace{0.15cm}  &~70\,500~   & ~--~  & ~--~ &~0.027~ & ~--~& ~--~&~--~ \\[4mm]
 		1D-ductile\hspace{0.15cm}  &~70\,500~   & ~--~  & ~--~ &~--~ & ~30~& ~250~& ~330~\\[4mm]
 		2D-brittle\hspace{0.15cm}  &~--~   & ~75\,100~  & ~28\,010~ &~25~ & ~--~ &~--~ & ~--~  \\[4mm]
 		2D-ductile\hspace{0.15cm}  &~--~   & ~136\,500~  & ~70\,300~ &~--~ & ~25~& ~300~&~443~\\[4mm]
 		3D-brittle\hspace{0.15cm}  &~--~    & ~121\,150~  & ~80\,770~ &~0.0027~ & ~--~&~--~& \\
 		\hline
 	\end{tabular}	\label{table1}
 \end{table}

\sectpa[Section4]{Numerical examples}

\begin{figure}
	\centering
	{\includegraphics[width=0.55\textwidth]{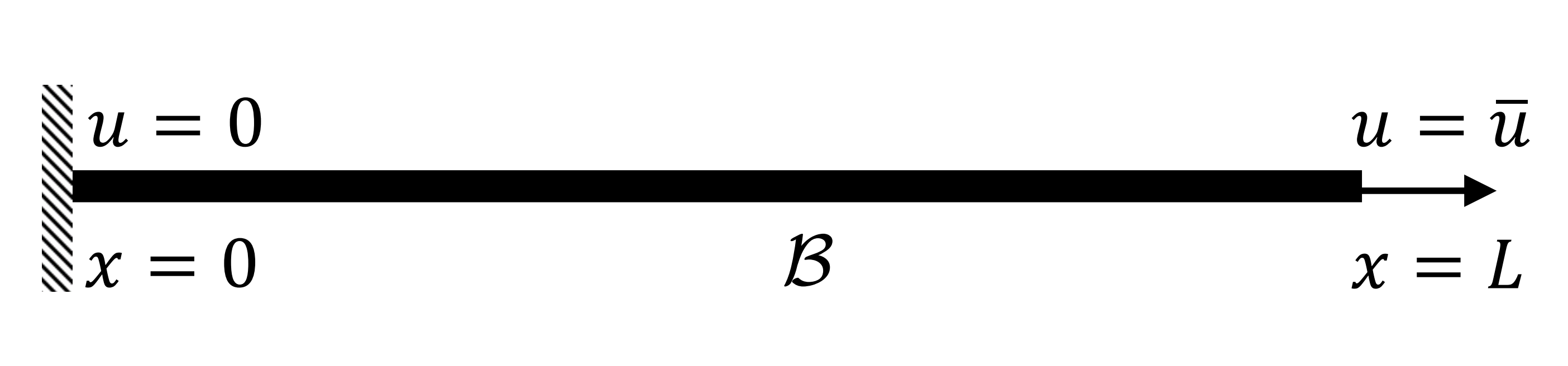}}  
	\caption{Example 1. Geometry and loading setup for the one-dimensional bar.}
	\label{BVP1D}
\end{figure}

This section demonstrates the performance of the proposed stochastic phase-field approach for brittle and ductile fracture in one-, two- and three-dimensional cases. As outlined in Section \ref{Section1}, the materials in standard phase-field problems are assumed to have a uniform macroscopic structure and properties (globally). However, those properties vary spatially at the heterogeneous microstructure (locally). In the local approach we monitor the propagation of crack through the matrix material. Whereas, in the global approach,  apart from the crack geometry, we will compute the expected value and the variance of the force response normal to the top boundary. For this, at time step $t\in\mathbb{T}$, we define the quantity $\BfrakJ$  as
\begin{align}
	\BfrakJ^{(i)}:=   	\BfrakJ(\Bu^{(i)},d^{(i)},\alpha^{(i)}) := \int_{\partial_D\calB\;|\;\overline{\Bu}\neq0} \bm{n} \cdot\bm{\sigma}^{(i)} \cdot \bm{n} \;\mathrm{d}{\bm{x}} \;\;\forall\;i\in N,
\end{align}
  where $\Bn$ is the outward unit vector normal to the Dirichlet boundary $\partial_D\calB$.

\begin{figure}[!b]
\includegraphics[width=\textwidth]{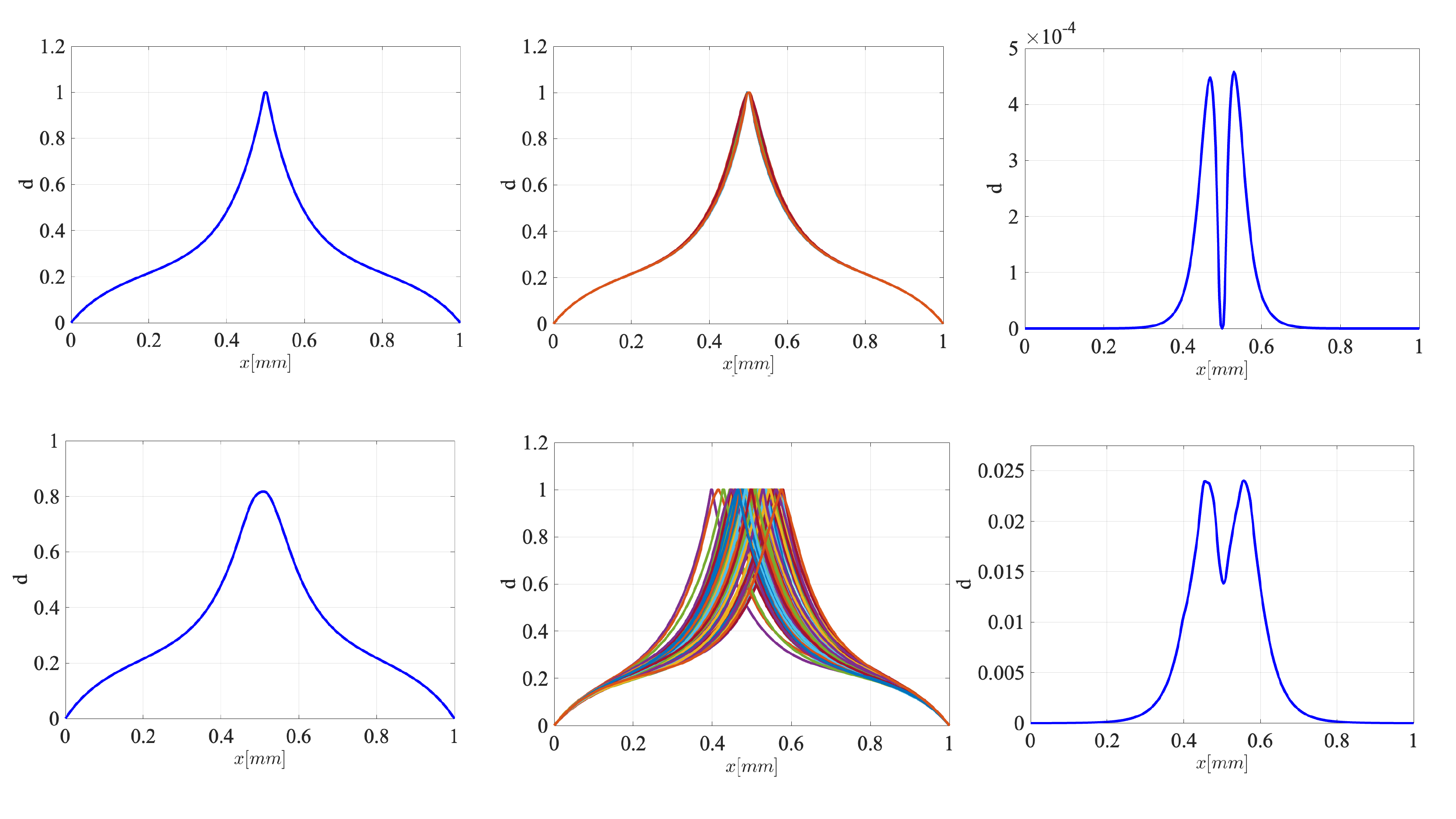}	
\caption{Example 1a (brittle fracture) with 10\% variation. The mean value (first column), {100 different samples} (middle), and the variance (third column) of the homogeneous case (first row) and the heterogeneous case (second row).}
\label{oneD_brittle}
\end{figure}

\sectpb[Section41]{One-dimensional brittle and ductile fracture}
 As the first case study, we consider the material (microscopic) fluctuation in brittle and ductile fracture. Hereby, a bar of unitary length $L = 1$ is considered where $x\in\calB:=[0,1]$ that is initially unstretched and undamaged. Its left end is fixed, i.e., $x=0$, while on its right end, i.e., $x=1$, a monotonic displacement increment  $\Delta \bar{u}=1\times10^{-5}$ mm is applied for 151 time steps. The example setup and boundary conditions are shown in Figure \ref{BVP1D}. Regarding to the finite element mesh size, 300 elements are used and {we replicate} the sampling in 400 iterations. The material parameters denoted as $\Bfrakq$ based on deterministic values $\bar{\Bfrakq}$ are given in Table \ref{table1}. 
%

For the stochastic case, we define two scenarios. First, materials have a homogeneous structure with a given variation. In each simulation, a value of the material parameter is determined using \eqref{theta}. The second possibility is related to the heterogeneity, in which at each point of the bar, the material parameters are fluctuated (i.e., the point-wise variation). In other words, we use \eqref{theta} to estimate the material values at each point.

 	\begin{figure}[!t]
\includegraphics[width=\textwidth]{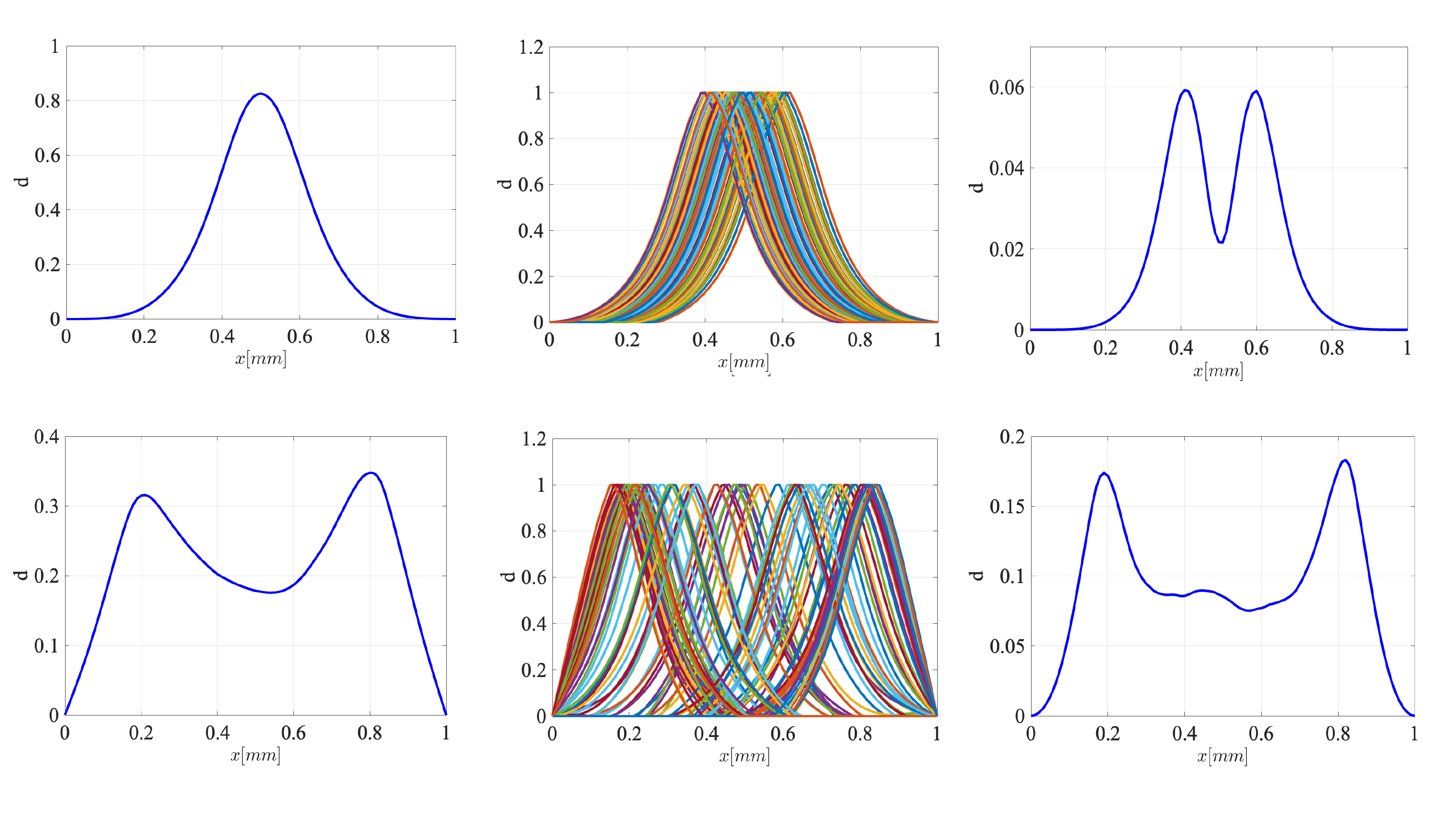}	 		
 		\caption{Example 1b (ductile Fracture) with 5\% variation. The mean value (first column), 100 different samples (middle), and the variance (third column) of homogeneous case (first row) and heterogeneous case (second row).}
 		\label{oneD_ductile1}
 	\end{figure}

\sectpc[Section4c]{Brittle fracture}

In this numerical test, the $\texttt{E}$lastic-$\texttt{D}$amage behavior in one-dimensional setting will be considered. Specifically, we investigate how the uncertainty affects the crack-surface. Figure \ref{oneD_brittle} (the first row) shows the mean value and $100$ different crack patterns (using $400$ simulations) for the homogeneous case. As the crack phase-field profile is regularized (using the length scale), only a minor variation around the peak point (at $x=0.5$ where the crack starts) is occurred. However, the fluctuation does not affect the fracture point (zero variance at $x=0.5$). In heterogeneous cases, the crack point has been varied due to the heterogeneity, and the variance is significantly higher compared to the homogeneous case.

\sectpc[Section4d]{Ductile fracture}
	
In this numerical example, we perform a stochastic analysis on the $\texttt{E}$lastic-$\texttt{P}$lastic-$\texttt{D}$amage behavior in one-dimensional setting of ductile fracture. Similar to the brittle fracture, two cases are considered, namely homogeneous and heterogeneous cases with a 5\% variation. The results are demonstrated in Figure \ref{oneD_ductile1}. Herein, although the crack profile is regularized, the imposed variation affects the crack-pattern significantly. However, the crack point shows a fixed value (negligible variation in the crack point). 
In contrast to the homogeneous case (first row), a negligible variation for the heterogeneous scenario will change the crack profile considerably. As Figure \ref{oneD_ductile1}(second row) shows, a significant variance is observed in the crack profile indicating the effect of the heterogeneity.
 	
%
%

 \sectpb[Section4_2D]{Two-dimensional microstructural RVE under tension}

\begin{figure}[!t]
	\centering
	{\includegraphics[clip,trim=2cm 3cm 0.4cm 5cm, width=14cm]{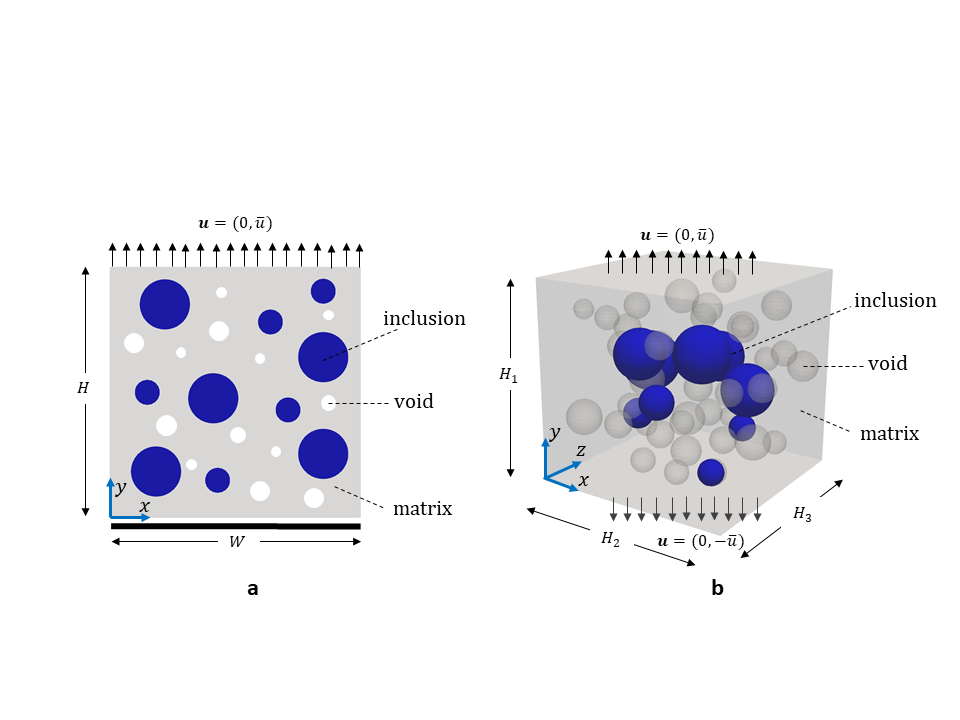}}  
	\caption{Geometry and loading setup for (a) the two-dimensional setting in Example 2, and (b) the three-dimensional setting in Example 3.}
	\label{BVP_23D}
\end{figure}

In the second case study, a two-dimensional microstructure RVE with stochastically distributed inclusions (aggregates) and voids (pores) under tension is considered. Both brittle ($\texttt{E-D}$), and ductile fractures ($\texttt{E-P-D}$) settings are considered. A boundary value problem applied to the square plate is shown in Figure \ref{BVP_23D}a. We set $H=1\;mm$ and $W=H$ hence $\calB=(0,1)^2$ $mm^2$ that includes randomly allocated inclusions and voids in domain. As a loading setup, we set the initial values for displacement and phase-field as $\bm u_0:=0 \in\calB$ and $d_0:=0 \in \calB$. Here, Galerkin finite element method with $H^1$-conforming bilinear (2D) elements are used for the $Q_1$-finite elements. A minimum element size of $h=0.01\;mm$ is considered such that the spatial discretization of the model includes approximately 20,000 four-node quadrilateral elements. Thus, the  fracture length-scale is set as $l=0.02\;mm$. The condition fulfills the heuristic requirement $h<l/2$ for the element size inside the localization zone (i.e., the support) of $d$, see \cite{miehe+welschinger+hofacker10}. Note that the plane-strain situation is considered. The displacement control is used with increments of $\Delta \bar{u}=1\times10^{-4}$ and 600 time steps. The material parameters are given in Table \ref{table1}.

To deal with heterogeneous microstructure materials, we further define the mismatched ratio between two categories of materials (inclusions and voids) pointing out by $\chi$, described as
 $$\bullet_\text{inc}=\chi\bullet_\text{mat}\quad \text{for}\quad \bullet\in\{E,\mu,K,G_c,\psi_c,\sigma_Y\}.$$
  In the upcoming examples, we set $\chi=10$, i.e., inclusions (aggregates) are $10$ times stiffer compared to the listed homogeneous structure in Table \ref{table1}. Next, we employ our allocating strategy to simulate the random distribution of voids and inclusions. This randomness procedure includes:
	\begin{enumerate}
	\item Random density of inclusions/voids: 
	\begin{itemize}
		\item Inclusions varies between 30 and 40 percent of the whole volume. 
		\item Voids varies between 5 and 10 percent of the whole volume. 
	\end{itemize}
	\item Random size of the particles (inclusions/voids).
	\begin{itemize}
		\item Inclusions varies between 5 mm and 20 mm (radius). 
		\item Voids varies between 0.2 and 5 mm (radius). 
	\end{itemize} 
	\item  Random position of the particles (see Algorithm \ref{Algorithm1}).
\end{enumerate} 
\begin{figure}
	\centering
	\subfloat{\includegraphics[width=0.65\textwidth]{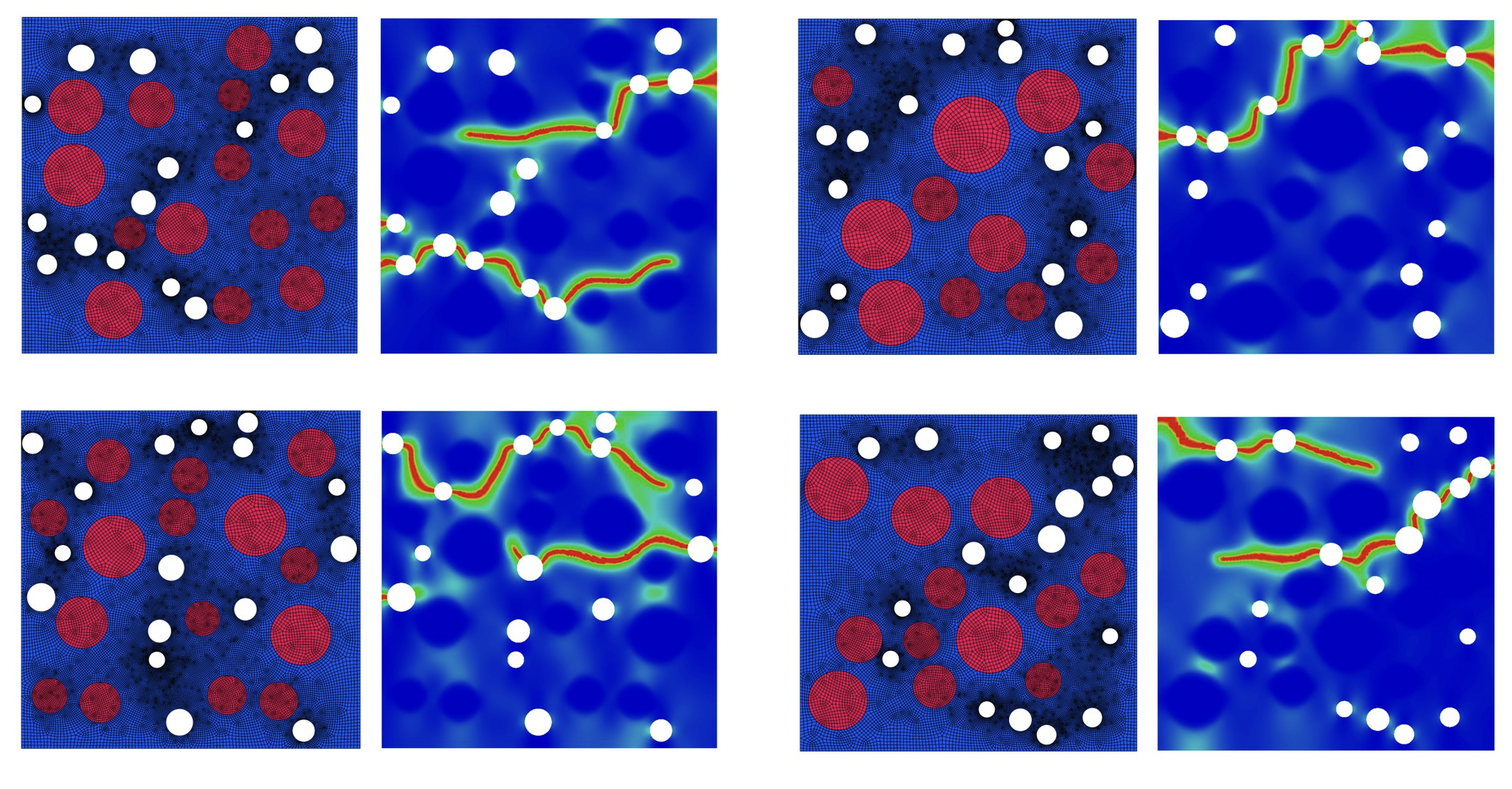}}\\[-0.5mm]
	\subfloat{\includegraphics[width=0.647\textwidth]{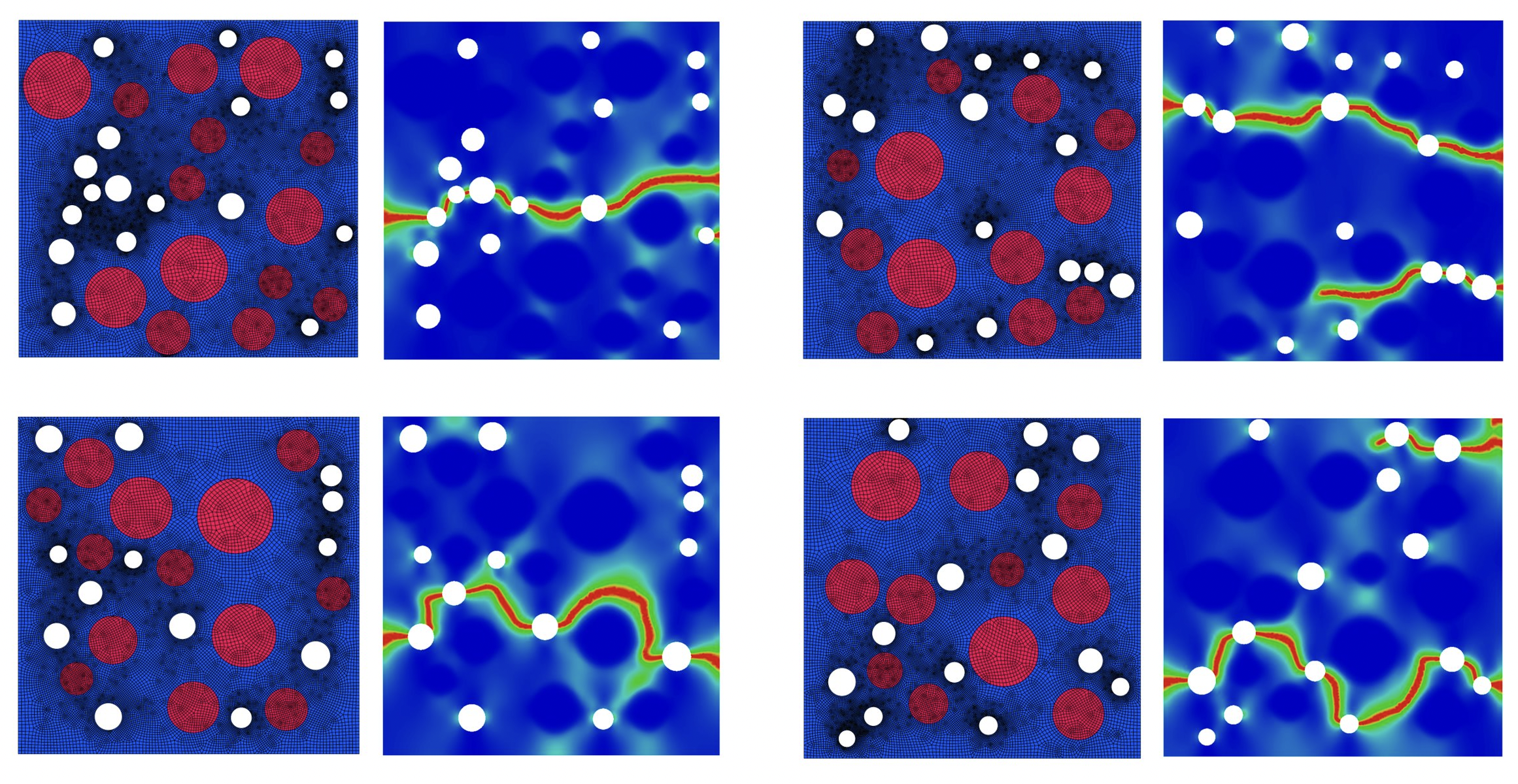}}\\[-0.5mm]
	\subfloat{\includegraphics[width=0.653\textwidth]{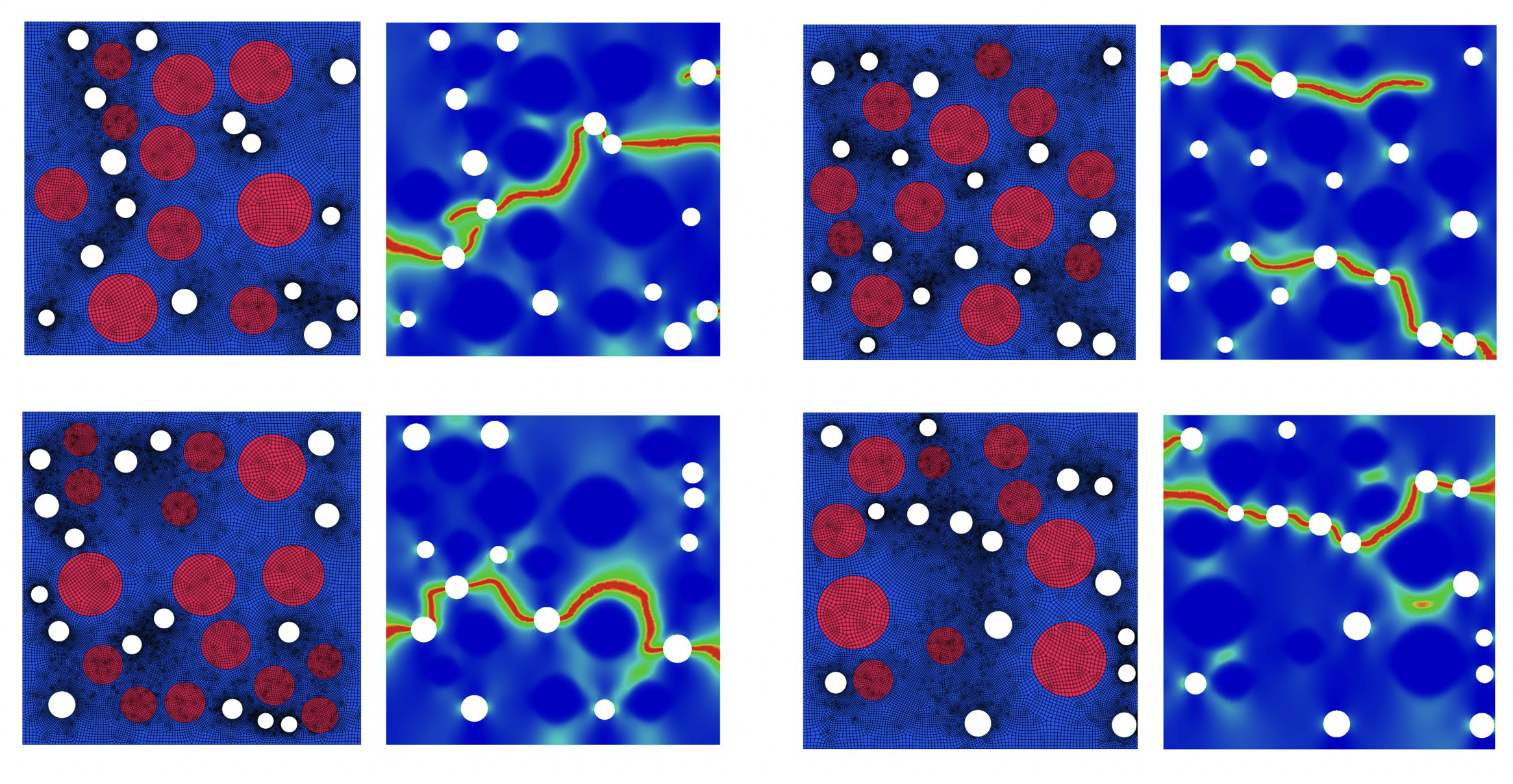}}\\[-0.5mm]
	\subfloat{\includegraphics[width=0.652\textwidth]{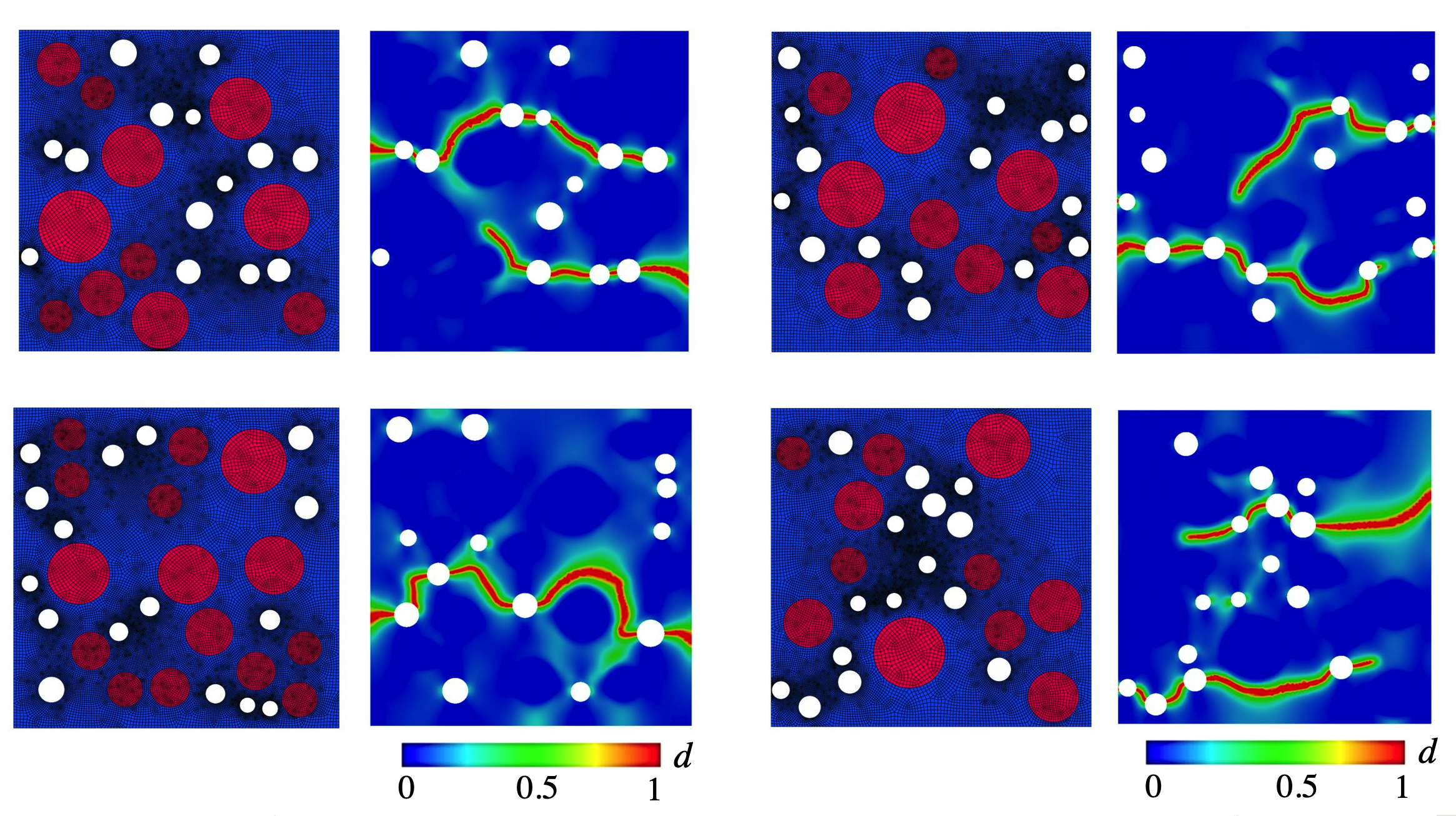}}
	\caption{Example 2. $16$ different mesh configuration denoting the random distribution of voids and inclusions in 2D and the corresponding crack pattern in brittle fracture.
	}
	\label{brittle_mesh1}
\end{figure}
\begin{figure}
	\centering
	\subfloat{{\includegraphics[clip,trim=1cm 3.2cm 1cm 0.5cm, width=15.5cm]{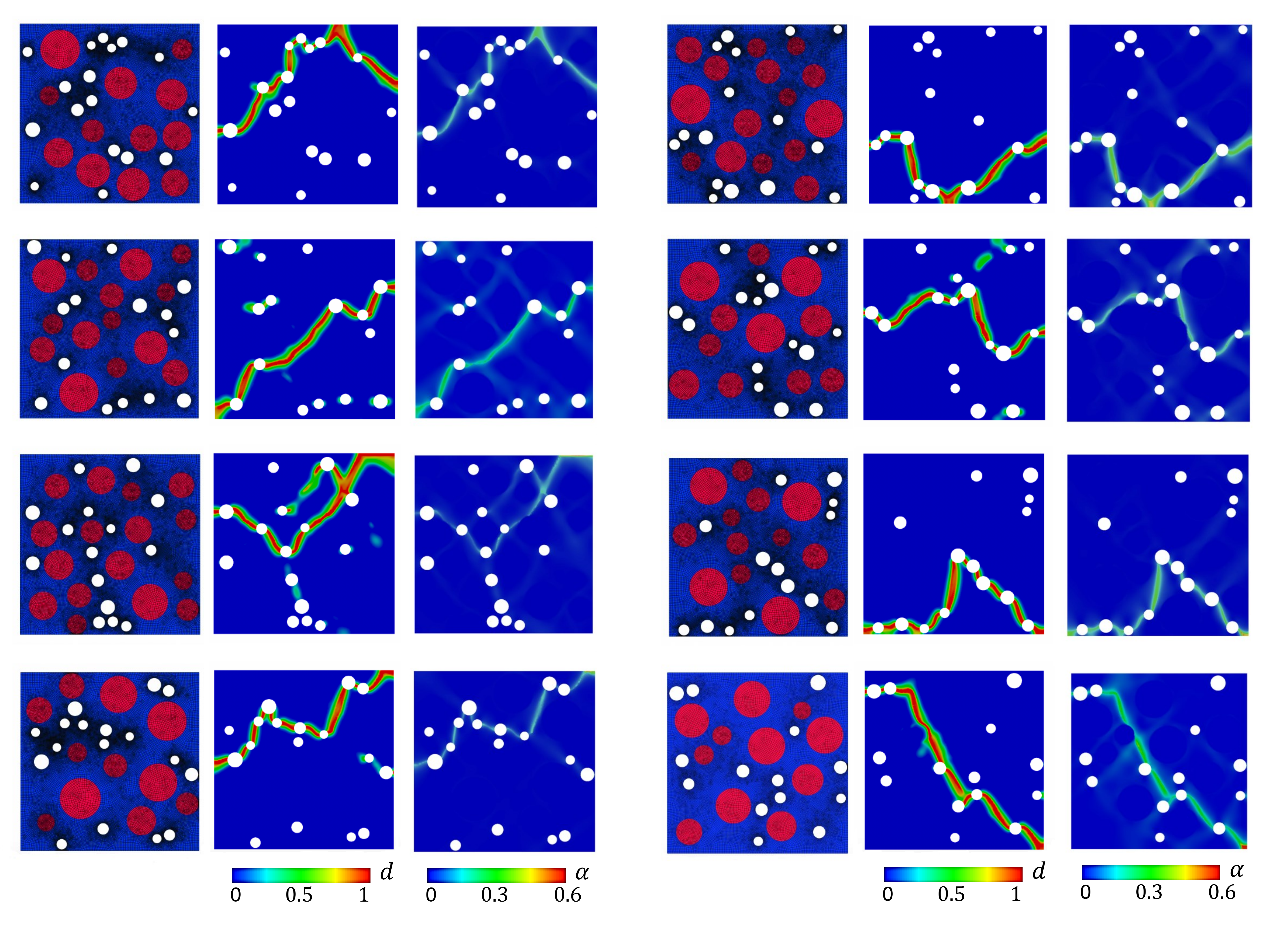}}}\\
	\subfloat{{\includegraphics[clip,trim=1cm 1cm 1cm 0.5cm, width=15.5cm]{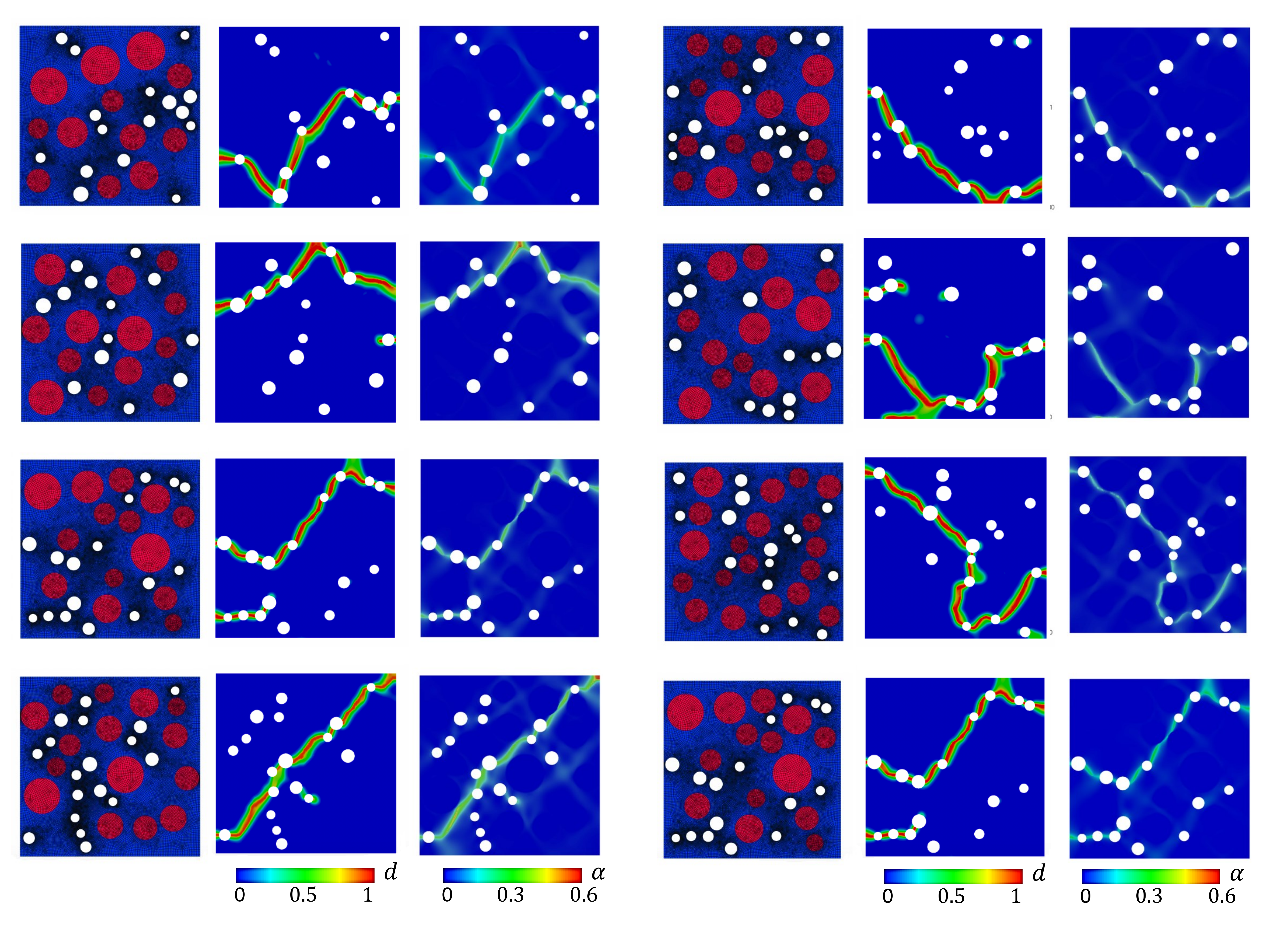}}} 
	\caption{Example 2. $16$ different mesh configuration denoting the random distribution of voids and inclusions in 2D along with the corresponding crack pattern and hardening in ductile fracture case.
	}
	\label{ductile_mesh1}
\end{figure}
\begin{figure}[!t]
	\subfloat{\includegraphics[width=0.5\textwidth]{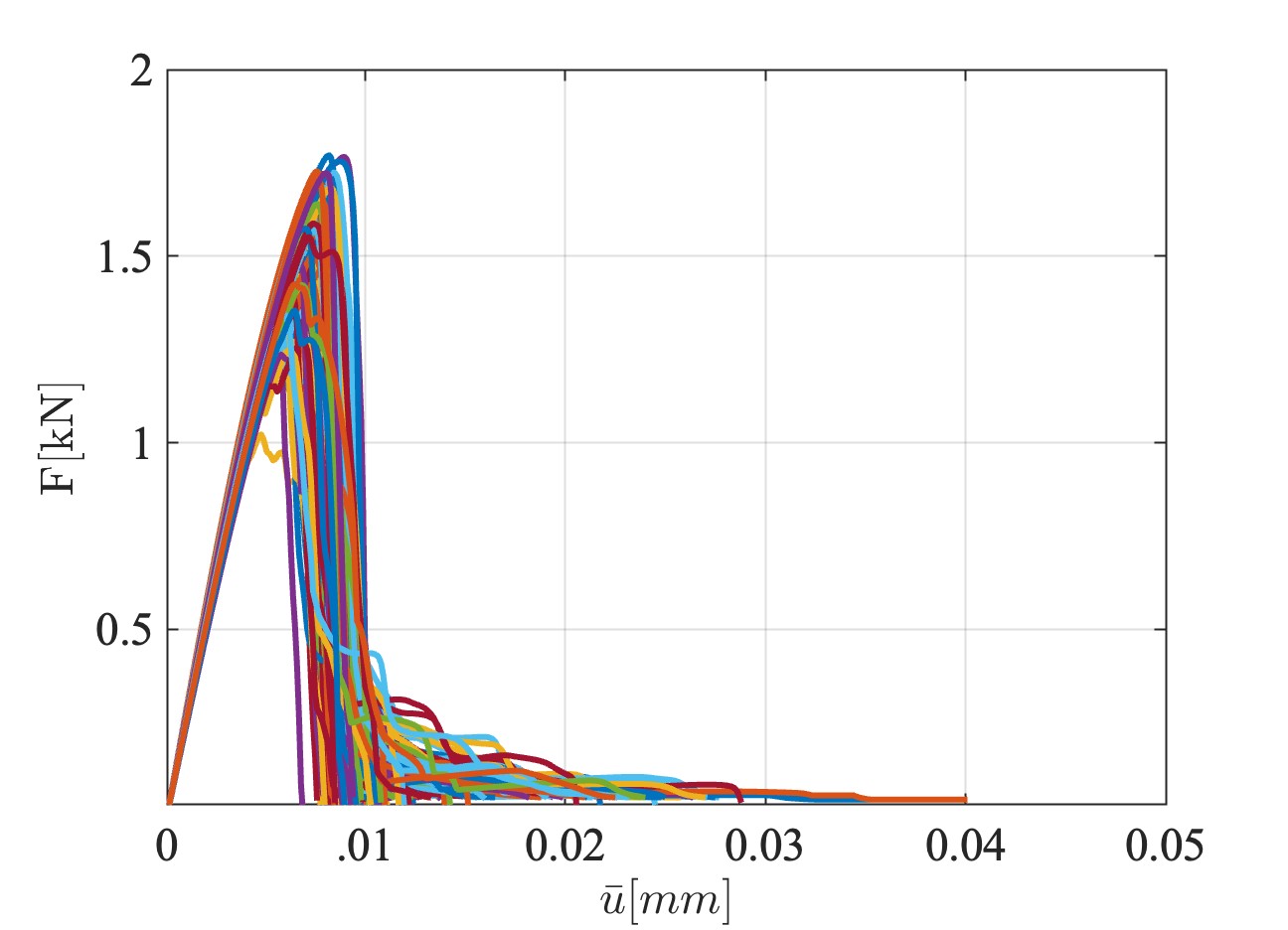}}
	\subfloat{\includegraphics[width=0.5\textwidth]{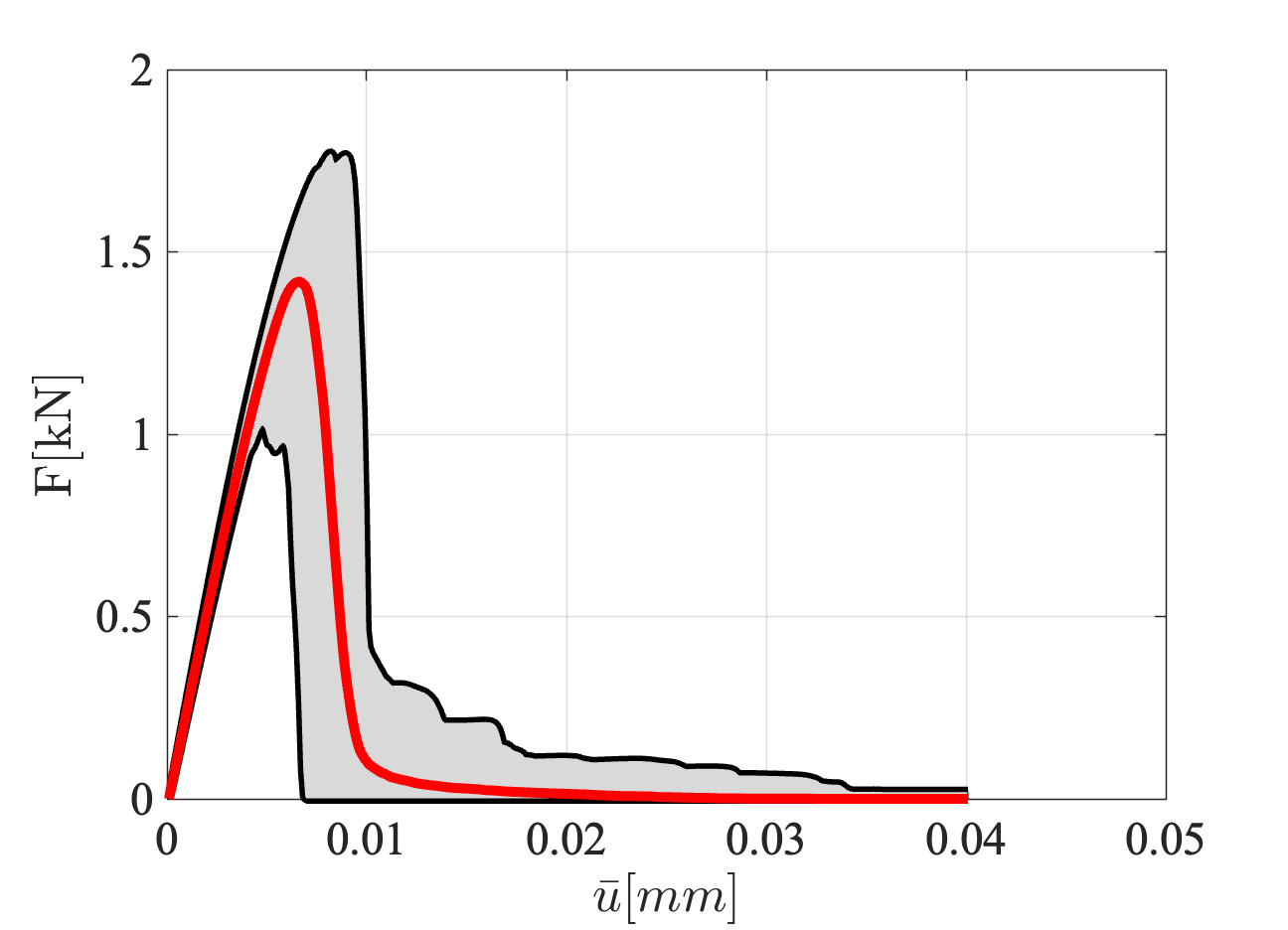}}\\
	\subfloat{\includegraphics[width=0.5\textwidth]{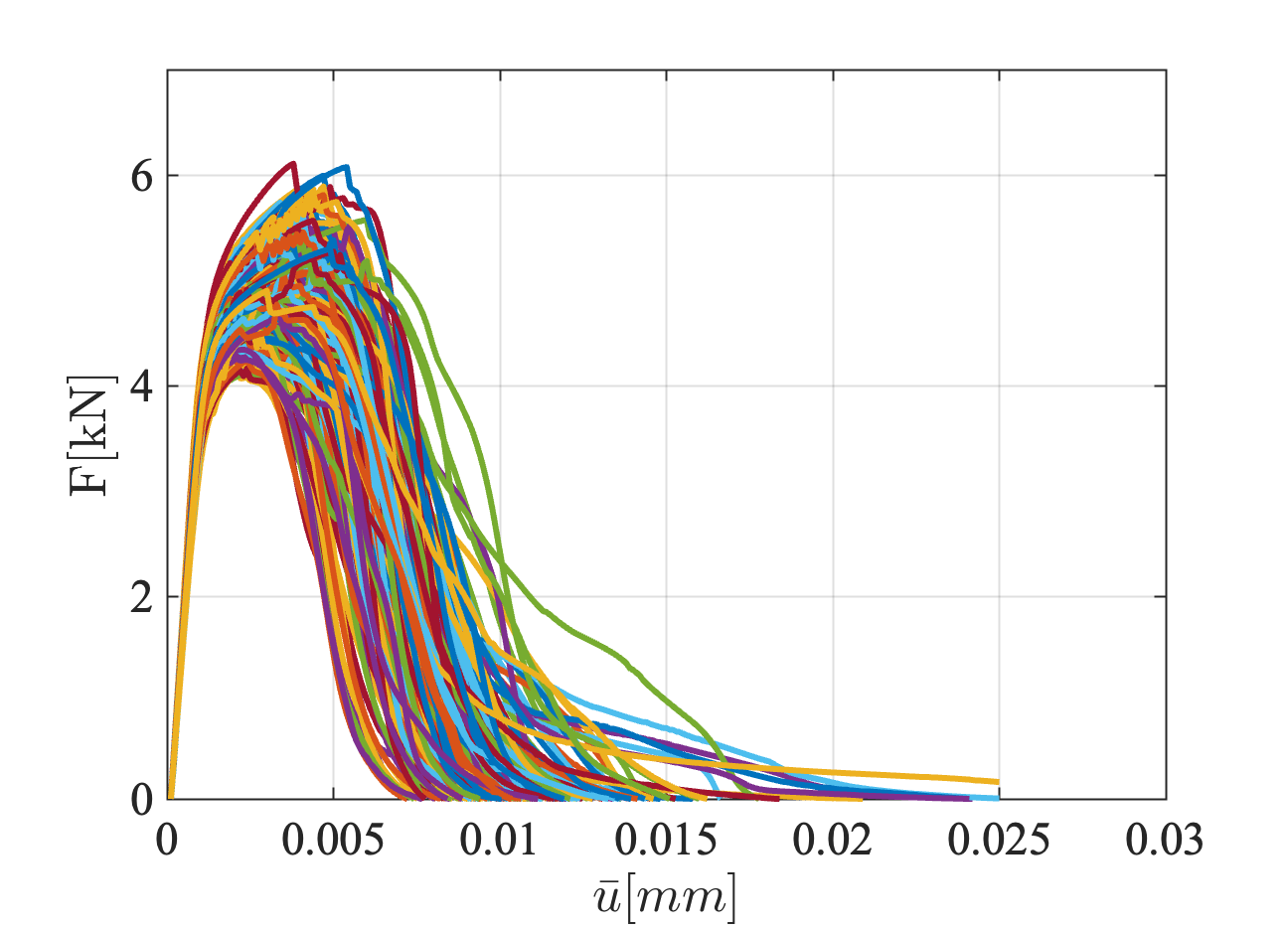}}
	\subfloat{\includegraphics[width=0.5\textwidth]{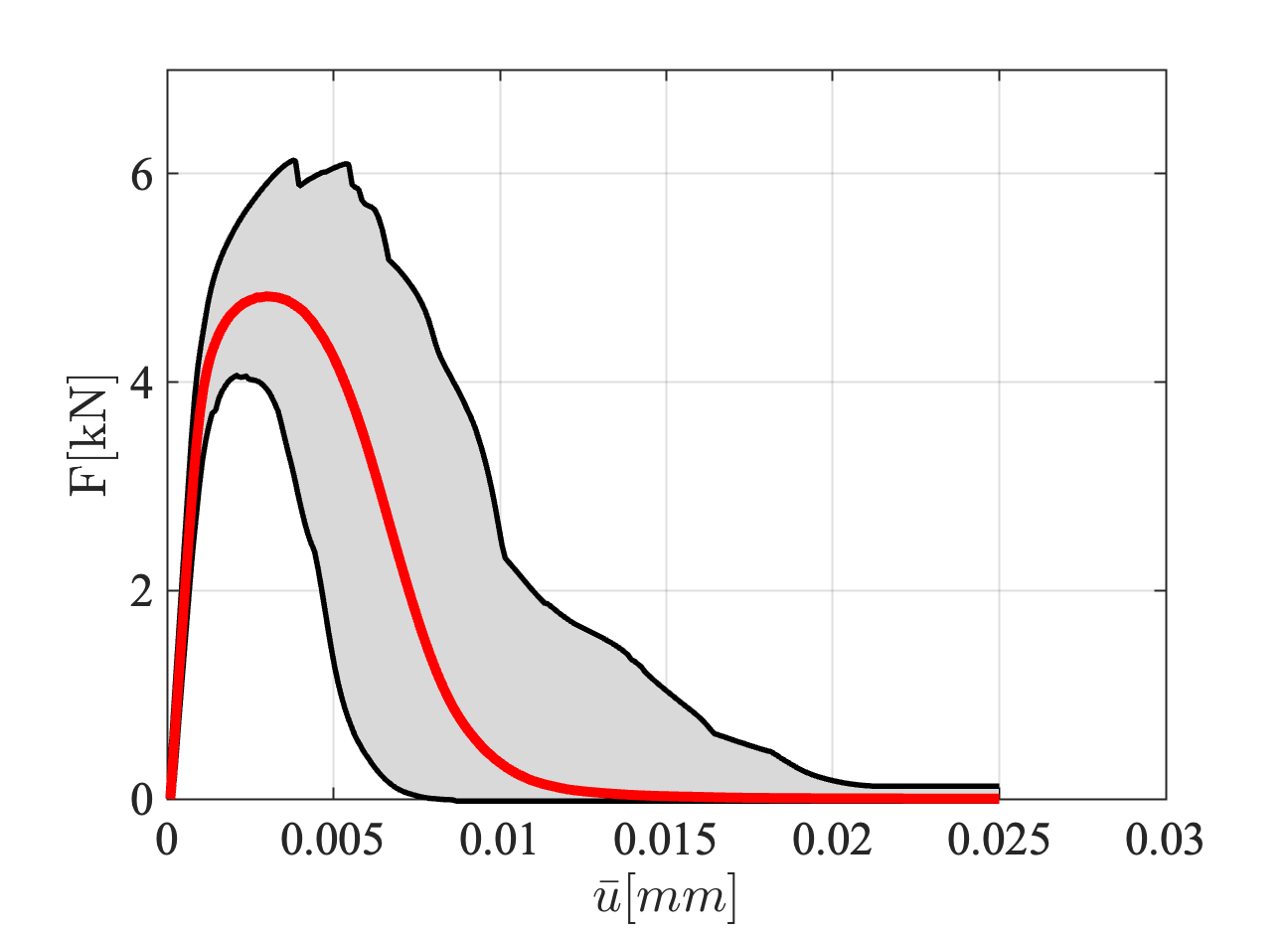}}
	\caption{Example 2. The load-displacement curves for 300 different random distribution of inclusions and voids for the brittle (top) and ductile (bottom) materials. The gray region (in the right column) shows the area between the maximum and the minimum of the diagrams during different time-steps. The mean values are shown with a solid red line.}
	\label{load2D}
\end{figure}
The random position of the particles will give rise to different crack propagation behavior. For instance, congestion of the inclusions (aggregates) in a part prevents the crack extension in this region. On the other hand, several voids will facilitate crack propagation easily. Figure  \ref{brittle_mesh1} shows the different distribution of the inclusions and voids. Considering the brittle case, we have the corresponding crack behavior in Figure \ref{brittle_mesh1}. As observed, the randomness in the matrix-material (e.g. concrete) leads to completely different crack patterns. In fact, after the crack nucleation, it propagates through the voids and among the inclusions. Considering the ductile concrete, $16$ different mesh configurations along with the related fracture behavior and hardening are illustrated in Figure \ref{ductile_mesh1}. 

In order to study the crack behavior during different time-steps, the load-displacement curve for 300 different random distributions are plotted for both brittle and ductile structure materials. Furthermore, we define a region denoting the maximum and minimum of the load-displacement curves for these simulations. The obtained information shows the possible range for all events. Figure \ref{load2D} demonstrates the diagrams, ranges and the mean values (shown in red). As shown, the ductile materials are much more resistant to fracture compared to the brittle materials, i.e. more than $3$ times of the force is needed. Furthermore, in the brittle case, the fracture happens sharply; whereas, in the ductile case, the crack requires more time to initiate (due to the plastic deformation).

%

\begin{table}
	\caption{The three-dimensional example: Two defined cases considering the sources of uncertainty in inclusions/voids and material parameters. In \textbf{Case a}, the heterogeneous structure using a variation $\eta=10\%$ in the material parameters is considered. In \textbf{Case b}, the first row is related to inclusions and the second one is for voids. Here, $\calU$ denotes the uniform distribution.}
	\vspace{1mm}
	\centering
	\begin{tabular}{|l |c c c c  |}
		\hline
		Uncertainty \hspace{0.3cm}&$~~\text{radius}~~$        & $~~\text{position}~~$    &   $~~\text{density}~~$   &   $~~\text{materials}~~$         \\\hline
		\textbf{Case a} \hspace{0.1cm}  &  $ \text{constant}$  & \text{constant} &  \text{constant} &$\eta=10\%$   \\[2mm]	\multirow{2}{4em}{\textbf{Case b } }
		\hspace{0.1cm}  &$\calU(6,15)$    & \text{random}    & $\calU(10,25)$  &$\text{constant}$    \\[2mm]
		\hspace{0.1cm}  &$\calU(4,8)$    & \text{random}    & $\calU(5,12)$  &$\text{constant}$    \\
		\hline
	\end{tabular} \label{tableAB}
\end{table}

 \sectpb[Section4_3D]{Three-dimensional microstructural RVE under tension}
 
 \begin{figure}[!b]
 	\subfloat{\includegraphics[width=1\textwidth]{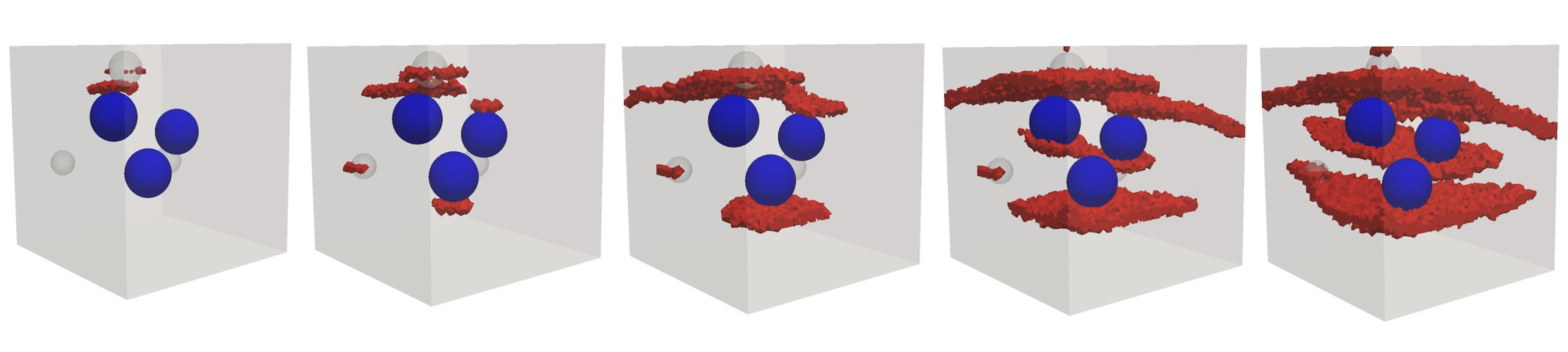}}\\ 
 	\caption{Example 3 (\textbf{Case a}). The evolution of the crack phase-field for the 3D distribution of voids and inclusions in a heterogeneous case.
 	}
 	\label{crack3}
 \end{figure}
 
 \begin{figure}[!t]
 	\subfloat{\includegraphics[width=1.02\textwidth]{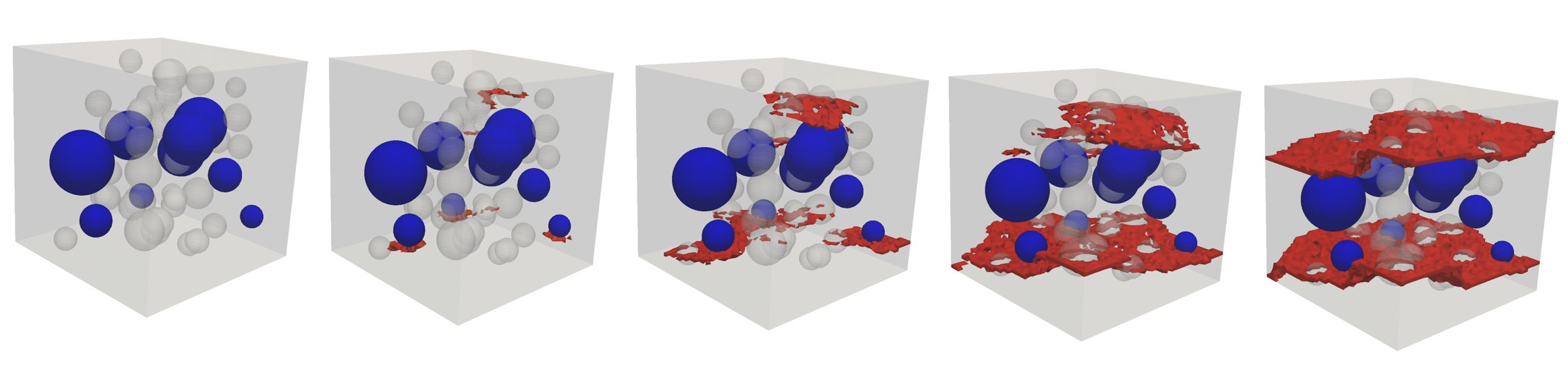}}\\ 
 	\caption{Example 3 (\textbf{Case b}). The evolution of the crack phase-field for the 3D distribution of voids and inclusions.
 	}
 	\label{stochastic1}
 \end{figure}
 
 In the last case study, a three-dimensional microstructure RVE with stochastically distributed aggregated and pores under tension is considered.  The three-dimensional setting helps us to monitor the crack propagation more illustratively. In this section, we only consider the brittle fracture ($\texttt{E-D}$) behavior, to avoid repetition compared to the analysis introduced in previous sections.
 
 A boundary value problem applied to the block specimen is shown in Figure \ref{BVP_23D}b. This is a tension test such that monomaniacal load is applied in both top and bottom (in opposite) directions. We set $H_1=1\;mm$ and $H_1=H_2=H_3$, and hence the cube space is $\calB=(0,1)^3$  that includes randomly allocated inclusions (aggregates) and voids (pores). As a loading setup, the initial values for displacement and phase-field are $\bm u_0:=0 \in\calB$ and $d_0:=0 \in \calB$. For the element technology, Galerkin finite element method with $H^1$-conforming trilinear (3D) elements is employed for FEM simulations. In this regard, the density of the inclusions varies between 10 and 25 percent, and the void density varies from 5 and 12 percent of the whole structure (concrete). The radii are from $6\,mm$ to $15\,mm$ for the inclusions and between $4\,mm$ and $8\,mm$ for the voids. Hereby, a tetrahedral meshes with an element size of $h=0.02\,mm$ is considered with averagely $500\,000$ elements in each simulation. A prescribed load of $\bar{u}=2\times 10^{-4}$ with 250 time-steps is used in the numerical simulation. The material properties are given in Table \ref{table1}, and $\chi=10$. In the following, we consider two different cases regarding the material fluctuation, described in Table \req{tableAB}.

 
%

\begin{figure} [!b]
	\begin{center}
		{\includegraphics[width=0.49\textwidth]{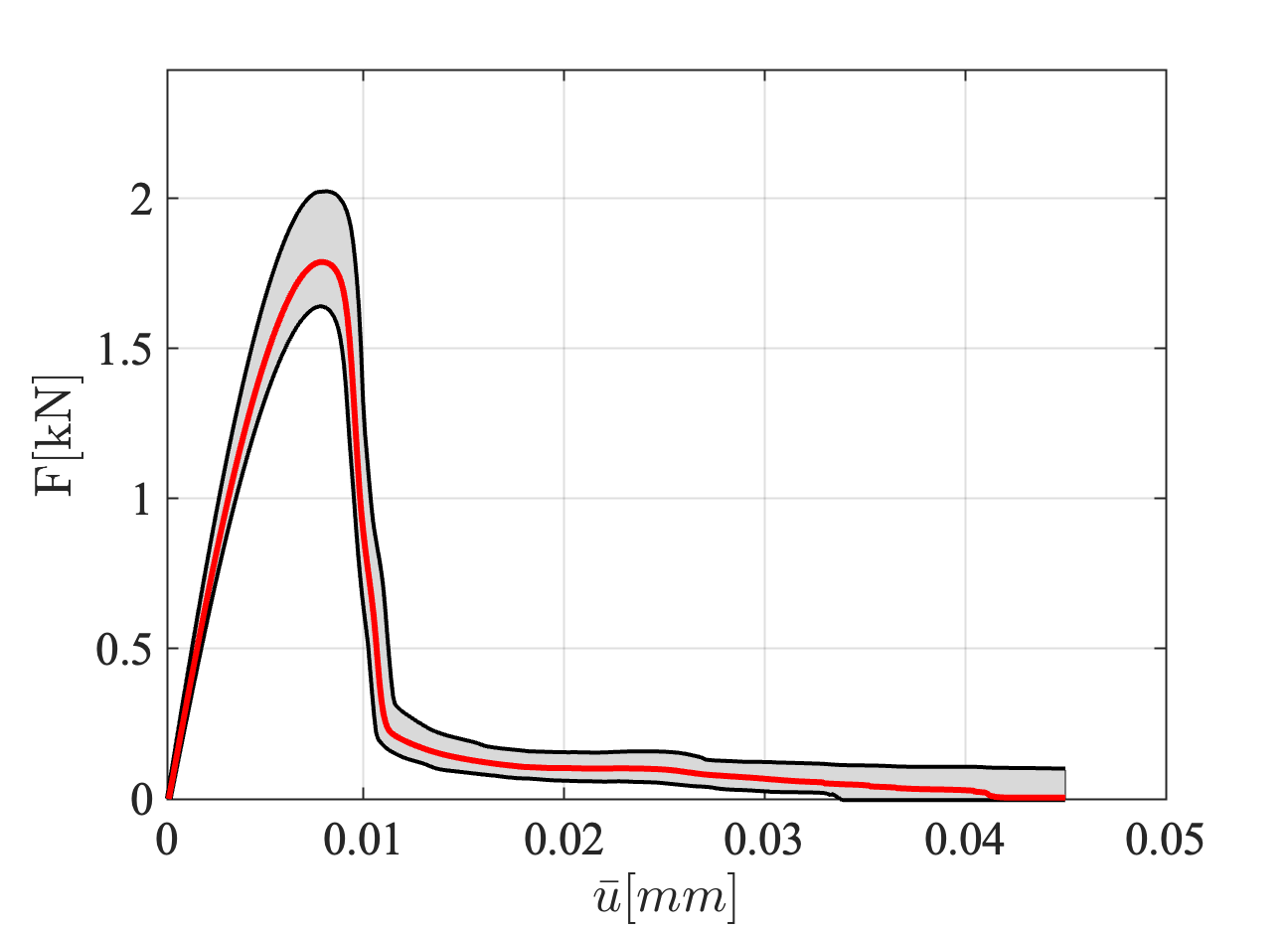}}  
		\includegraphics[width=0.49\textwidth]{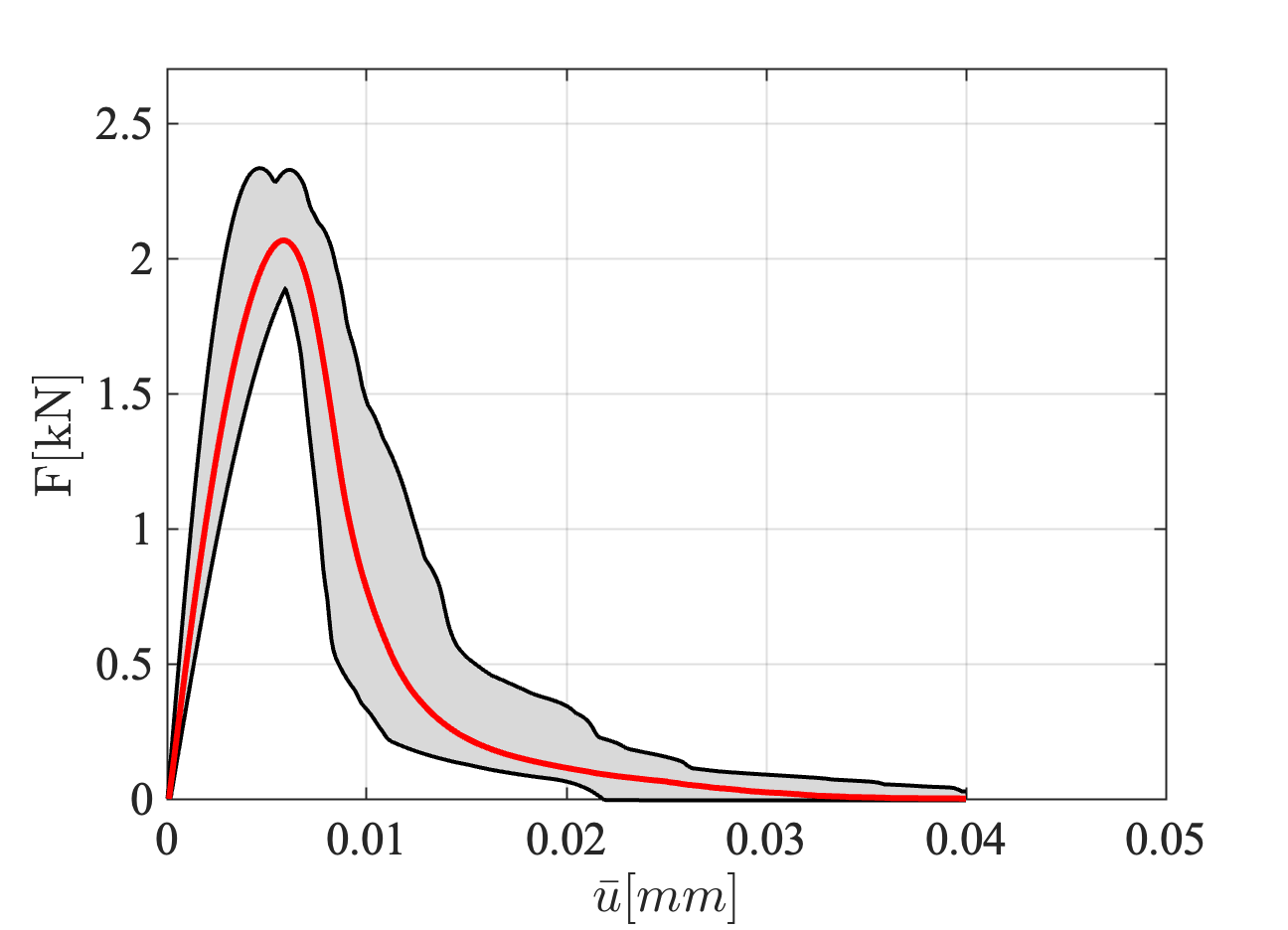}
	\end{center}
	\caption{Example 3. The load-displacement curve. \textbf{Case a} (left): Random heterogeneous material structure in a fixed spatial coordinates of different phases. \textbf{Case b} (right): Different random distribution of inclusions and voids using 300 replications for the brittle materials. The gray region shows the area between the maximum and the minimum of the diagrams during different time-steps. The mean values are plotted with a solid red-line.}
	\label{load3D}
\end{figure}

\begin{figure}[!t]
	\centering
	\subfloat{{\includegraphics[width=0.8\textwidth]{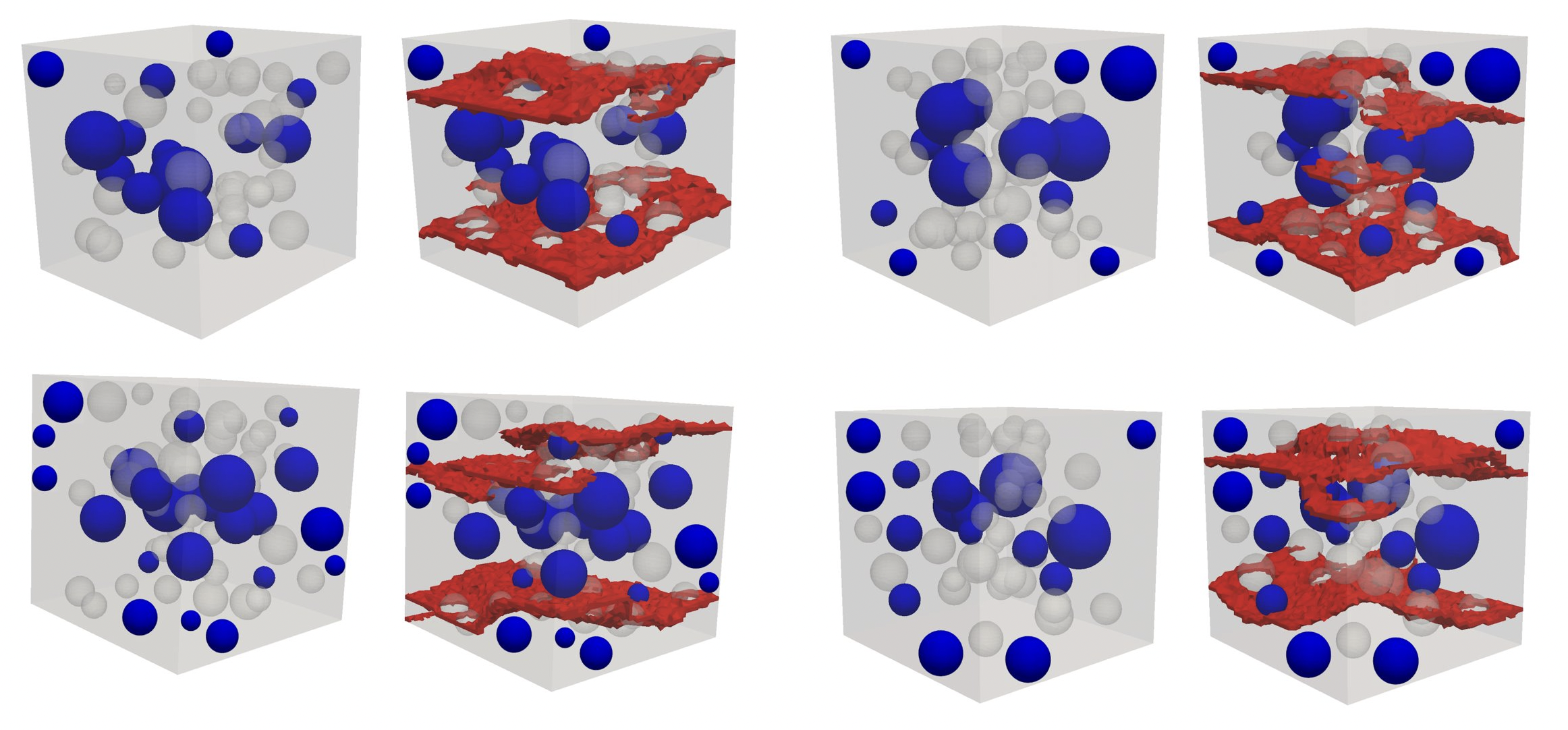}}}\\[-2mm]
	\subfloat{{\includegraphics[width=0.8\textwidth]{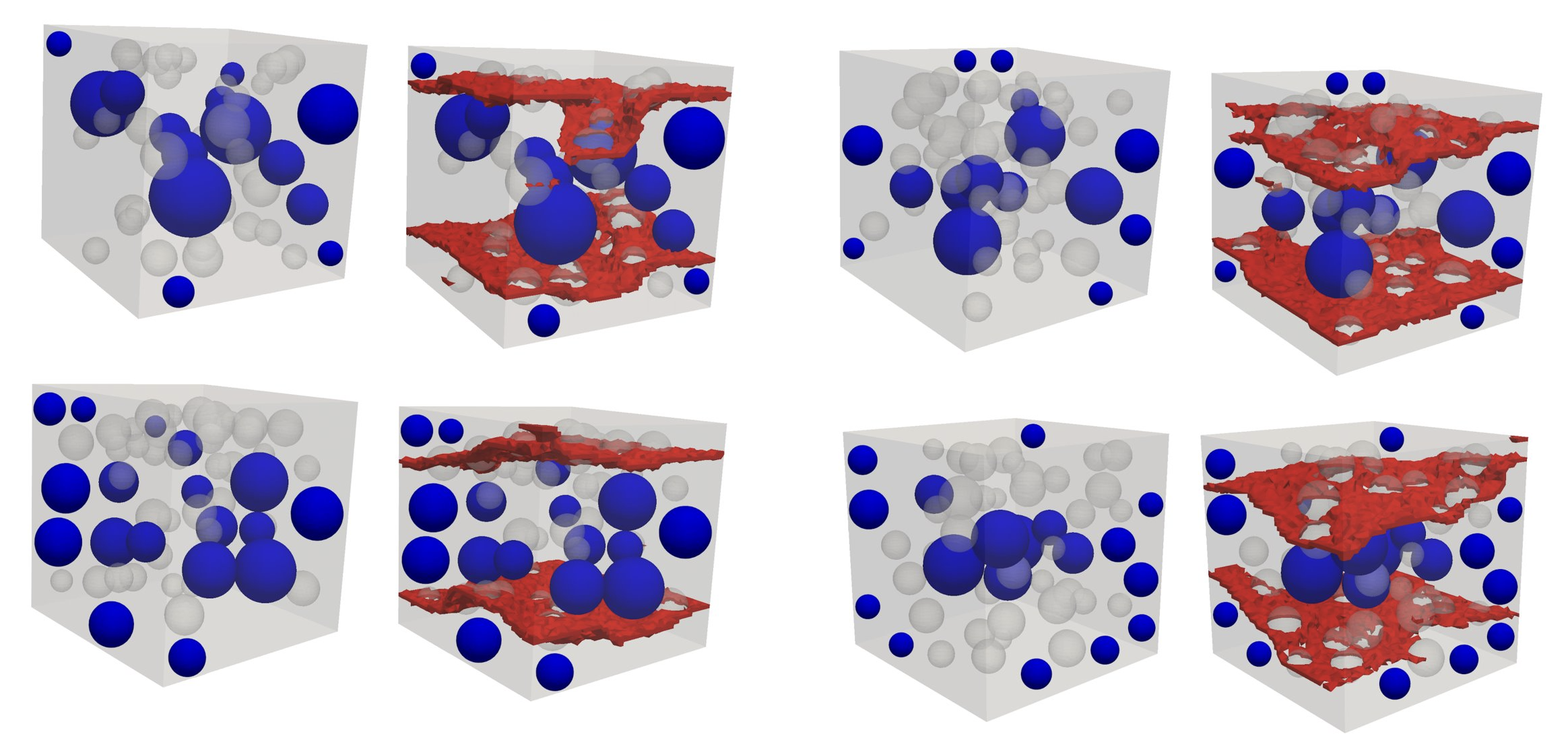}}}\\[-2mm]
	\subfloat{{\includegraphics[width=0.8\textwidth]{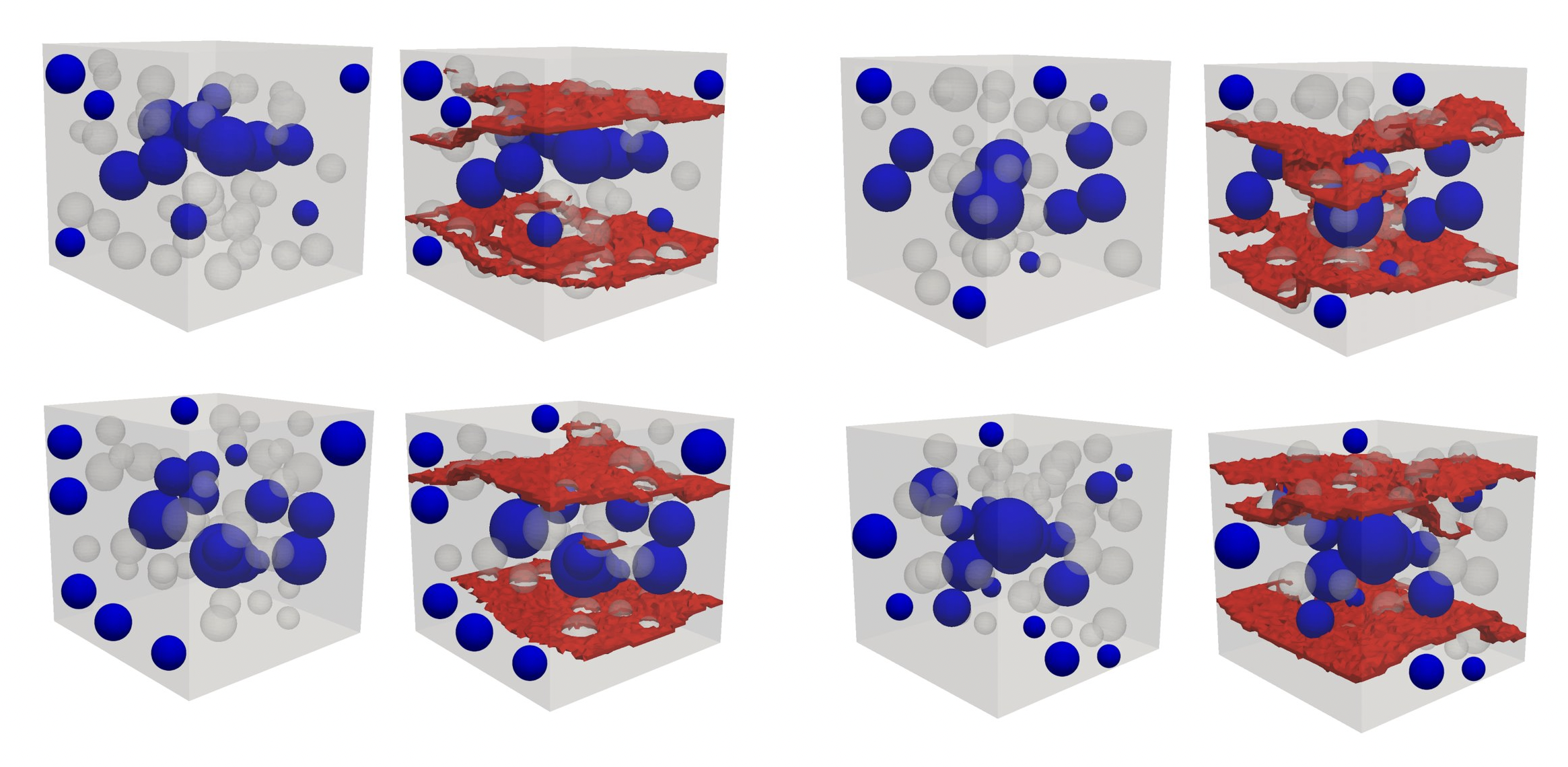}}}\\[-2mm]
	\subfloat{{\includegraphics[width=0.8\textwidth]{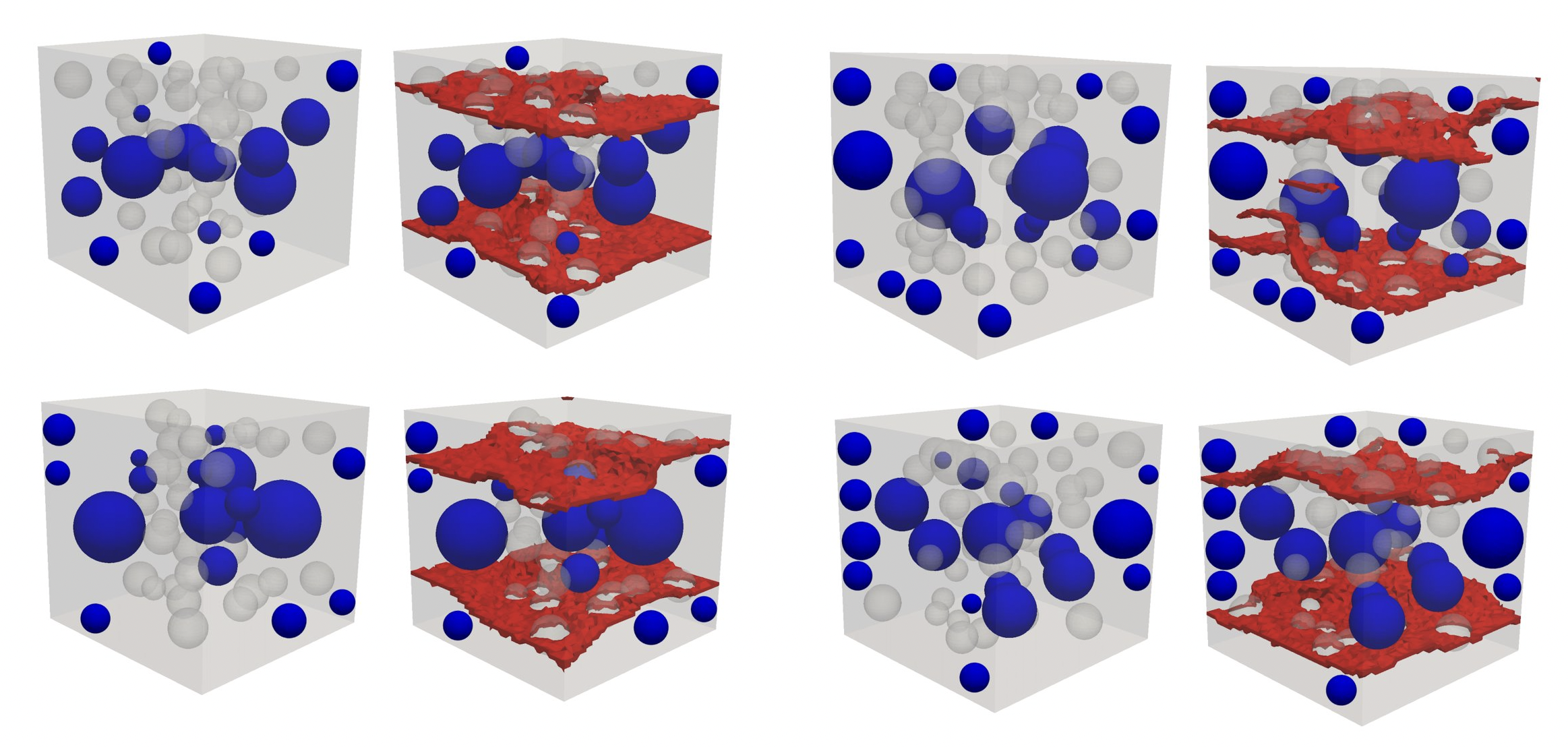}}}
	\caption{Example 3 (\textbf{Case b}). The crack phase-field for randomly distributed aggregates/pores in cube at the complete failure for the brittle fracture.}
	\label{3D2}
\end{figure}

\sectpc{Case a -  Random Material Properties}
First, we consider a heterogeneous structure using the material parameters variation of $\eta=10\%$ applied in \eqref{theta}. However, the inclusions/voids have no spatial fluctuation. Specifically, for different samplings, we have only a point-wise material variation with a constant position and number of the particles (inclusions/voids).  

The crack propagation during time is shown in Figure \ref{crack3} for an arbitrary sample. The variation interval for the load-displacement diagram using 300 simulations is depicted in \ref{load3D} (left). As only a variation in the material parameters is considered, the gray area of all possible solutions is smaller than other examples.

\sectpc{Case b - Geometrical Perturbation}
Next, different distribution of particles are investigated, in which the inclusions are $10$ times stiffer than the matrix material. We monitor the crack propagation during different time-steps, as shown in Figure \ref{stochastic1} for an arbitrary sample. The crack initiates in different parts of the specimen and propagates until the complete failure. In this analysis, the simulations are replicated 300 times. Figure \ref{load3D} (right) shows the maximum and minimum of the load due to uncertainties versus prescribed displacement and the corresponding expected value. 

To present the effect of uncertainties due to spatial variations of different phases, $16$ specific three-dimensional distributions of inclusions (aggregates) and voids (pores) are shown in Figure \ref{3D2} at the final deformation states. As expected, the randomness results in different fracture patterns.

\sectpa[Section5]{Conclusions}

 Heterogeneous materials at the lower scale are typically subjected to several uncertainties that inherently exist through the volume fraction defects at the micro or mesoscale. The classical approach to formulate those defects relies on a deterministic approximation of failure response, while such effects are not captured for the unavoidable uncertainties of each parameter associated with experimental observations. To overcome that, the current work is devoted to a rigorous mathematical formulation of stochastic-based variational formulations of failure mechanisms at the micro/meso-level. More specifically, uncertainties in brittle and ductile failure are investigated. The primary objective of this contribution is to model randomness and fluctuations of different phases in the highly heterogeneous meso/microstructures. 
 
To explore the fundamental nature of the proposed model, first, we studied a localization effect within a one-dimensional bar due to the variation in material proprieties for gradient-based plasticity and damage models. The main observation is that point-to-point correlations of
the crack phase fields in the underlying heterogeneous bar can be captured. These stochastic solutions are represented by random fields or random variables in contrast to the classical deterministic solution spaces. Next, in two- and three-dimensional scenarios,  by using the Monte Carlo finite element method, we modeled the random distribution of the inclusions/voids and considered their effects on the material stiffness locally (the crack propagation pattern in different slides) and globally (considering the force-displacement diagram).  In this way, different evolution of cracks at the lower scale emerges as a consequence of the underlying uncertainty of physical model parameters. To formulate these uncertainties, we developed a procedure for the allocating process of highly numbers of inclusions/voids with different volume fractions in such a way that there is no intersection between them. The results enable us to provide a confidence interval for the fracture energy denoting the minimum/maximum necessary force for the fracture. Hereby, the computed expected value represents the average for all heterogeneities.

 \subsection*{Acknowledgment}
Fadi Aldakheel (FA) appreciates the scientific support of the German Research Foundation in the Priority Program \texttt{SPP 2020} (project number: 353757395). FA would like to thank M.Sc. {\it Markus H\"upgen} and Prof. {\it Michael Haist} for providing the computer tomography CT-images of concrete mesostructure (Figure \ref{Figure1}), which represents an application of the segmentation tolerance (perturbations) in the geometrical properties.
 
 
%
	\clearpage
	 {\normalsize
	 	\begin{spacing}{0.8}
	 		\bibliographystyle{ieeetr}
	\bibliography{./lit}

\begin{thebibliography}{10}

\bibitem{zohdi2004introduction}
T.~I. Zohdi and P.~Wriggers, {\em An introduction to computational
  micromechanics}, vol.~20.
\newblock Springer Science \& Business Media, 2004.

\bibitem{wriggers2007micro}
P.~Wriggers and M.~Hain, ``Micro-meso-macro modelling of composite materials,''
  in {\em Computational Plasticity}, pp.~105--122, Springer, 2007.

\bibitem{hain2008numerical}
M.~Hain and P.~Wriggers, ``Numerical homogenization of hardened cement paste,''
  {\em Computational Mechanics}, vol.~42, no.~2, pp.~197--212, 2008.

\bibitem{wessels2022computational}
H.~Wessels, C.~B{\"o}hm, F.~Aldakheel, M.~H{\"u}pgen, M.~Haist, L.~Lohaus, and
  P.~Wriggers, ``Computational homogenization using convolutional neural
  networks,'' in {\em Current Trends and Open Problems in Computational
  Mechanics}, pp.~569--579, Springer, 2022.

\bibitem{BourFraMar08}
B.~Bourdin, G.~Francfort, and J.-J. Marigo, ``The variational approach to
  fracture,'' {\em Journal of Elasticity}, vol.~91, pp.~5--148, 2008.

\bibitem{KUHN10}
C.~Kuhn and R.~M\"uller, ``A continuum phase field model for fracture,'' {\em
  Engineering Fracture Mechanics}, vol.~77, no.~18, pp.~3625 -- 3634, 2010.

\bibitem{miehe+welschinger+hofacker10}
C.~Miehe, F.~Welschinger, and M.~Hofacker, ``Thermodynamically consistent
  phase-field models of fracture: Variational principles and multi-field {FE}
  implementations,'' {\em International Journal for Numerical Methods in
  Engineering}, vol.~83, pp.~1273--1311, 2010.

\bibitem{hesch+weinberg14}
C.~Hesch and K.~Weinberg, ``Thermodynamically consistent algorithms for a
  finite-deformation phase-field approach to fracture,'' {\em International
  Journal for Numerical Methods in Engineering}, vol.~99, pp.~906--924, 2014.

\bibitem{Wick15Adapt}
T.~Heister, M.~F. Wheeler, and T.~Wick, ``A primal-dual active set method and
  predictor-corrector mesh adaptivity for computing fracture propagation using
  a phase-field approach,'' {\em Computer Methods in Applied Mechanics and
  Engineering}, vol.~290, pp.~466 -- 495, 2015.

\bibitem{rezaei2021direction}
S.~Rezaei, J.~R. Mianroodi, T.~Brepols, and S.~Reese, ``Direction-dependent
  fracture in solids: Atomistically calibrated phase-field and cohesive zone
  model,'' {\em Journal of the Mechanics and Physics of Solids}, vol.~147,
  p.~104253, 2021.

\bibitem{alessi2020phase}
R.~Alessi, F.~Freddi, and L.~Mingazzi, ``Phase-field numerical strategies for
  deviatoric driven fractures,'' {\em Computer Methods in Applied Mechanics and
  Engineering}, vol.~359, p.~112651, 2020.

\bibitem{denli2020phase}
F.~A. Denli, O.~G{\"u}ltekin, G.~A. Holzapfel, and H.~Dal, ``A phase-field
  model for fracture of unidirectional fiber-reinforced polymer matrix
  composites,'' {\em Computational Mechanics}, vol.~65, no.~4, pp.~1149--1166,
  2020.

\bibitem{heider2021review}
Y.~Heider, ``A review on phase-field modeling of hydraulic fracturing,'' {\em
  Engineering Fracture Mechanics}, vol.~253, p.~107881, 2021.

\bibitem{schreiber2020phase}
C.~Schreiber, C.~Kuhn, R.~M{\"u}ller, and T.~Zohdi, ``A phase field modeling
  approach of cyclic fatigue crack growth,'' {\em International Journal of
  Fracture}, vol.~225, no.~1, pp.~89--100, 2020.

\bibitem{steinke2019phase}
C.~Steinke and M.~Kaliske, ``A phase-field crack model based on directional
  stress decomposition,'' {\em Computational Mechanics}, vol.~63, no.~5,
  pp.~1019--1046, 2019.

\bibitem{arash2021finite}
B.~Arash, W.~Exner, and R.~Rolfes, ``A finite deformation phase-field fracture
  model for the thermo-viscoelastic analysis of polymer nanocomposites,'' {\em
  Computer Methods in Applied Mechanics and Engineering}, vol.~381, p.~113821,
  2021.

\bibitem{fantoni2019phase}
F.~Fantoni, A.~Bacigalupo, M.~Paggi, and J.~Reinoso, ``A phase field approach
  for damage propagation in periodic microstructured materials,'' {\em
  International Journal of Fracture}, pp.~1--24, 2019.

\bibitem{aldakheel2021feed}
F.~Aldakheel, R.~Satari, and P.~Wriggers, ``Feed-forward neural networks for
  failure mechanics problems,'' {\em Applied Sciences}, vol.~11, no.~14,
  p.~6483, 2021.

\bibitem{Wick+2020}
T.~Wick, {\em Multiphysics Phase-Field Fracture: Modeling, Adaptive
  Discretizations, and Solvers}.
\newblock De Gruyter, 2020.

\bibitem{pillai2020anisotropic}
U.~Pillai, S.~P. Triantafyllou, Y.~Essa, and F.~M. de~la Escalera, ``An
  anisotropic cohesive phase field model for quasi-brittle fractures in thin
  fibre-reinforced composites,'' {\em Composite Structures}, vol.~252,
  p.~112635, 2020.

\bibitem{heider2022self}
Y.~Heider, F.~Bamer, F.~Ebrahem, and B.~Markert, ``Self-organized criticality
  in fracture models at different scales,'' {\em Examples and Counterexamples},
  vol.~2, p.~100054, 2022.

\bibitem{rezaei2022anisotropic}
S.~Rezaei, A.~Harandi, T.~Brepols, and S.~Reese, ``An anisotropic cohesive
  fracture model: advantages and limitations of length-scale insensitive
  phase-field damage models,'' {\em Engineering Fracture Mechanics}, p.~108177,
  2022.

\bibitem{selevs2021general}
K.~Sele{\v{s}}, F.~Aldakheel, Z.~Tonkovi{\'c}, J.~Sori{\'c}, and P.~Wriggers,
  ``A general phase-field model for fatigue failure in brittle and ductile
  solids,'' {\em Computational Mechanics}, vol.~67, no.~5, pp.~1431--1452,
  2021.

\bibitem{ambatiphase}
M.~Ambati, J.~Heinzmann, M.~Seiler, and M.~K{\"a}stner, ``Phase-field modelling
  of brittle fracture along the thickness direction of plates and shells,''
  {\em International Journal for Numerical Methods in Engineering}, 2022.
\newblock https://doi.org/10.1002/nme.7001.

\bibitem{zhuang2022phase}
X.~Zhuang, S.~Zhou, G.~Huynh, P.~Aerias, and T.~Rabczuk, ``Phase field
  modelling and computer implementation: A review,'' {\em Engineering Fracture
  Mechanics}, p.~108234, 2022.

\bibitem{seiler2021phase}
M.~Seiler, S.~Keller, N.~Kashaev, B.~Klusemann, and M.~K{\"a}stner,
  ``Phase-field modelling for fatigue crack growth under laser shock
  peening-induced residual stresses,'' {\em Archive of Applied Mechanics},
  vol.~91, no.~8, pp.~3709--3723, 2021.

\bibitem{tan2022phase}
W.~Tan and E.~Mart{\'\i}nez-Pa{\~n}eda, ``Phase field fracture predictions of
  microscopic bridging behaviour of composite materials,'' {\em Composite
  Structures}, vol.~286, p.~115242, 2022.

\bibitem{ambati+gerasimov+lorenzis15}
M.~Ambati, T.~Gerasimov, and L.~De~Lorenzis, ``Phase-field modeling of ductile
  fracture,'' {\em Computational Mechanics}, vol.~55, pp.~1017--1040, 2015.

\bibitem{aldakheel2020microscale}
F.~Aldakheel, ``A microscale model for concrete failure in poro-elasto-plastic
  media,'' {\em Theoretical and Applied Fracture Mechanics}, vol.~107,
  p.~102517, 2020.

\bibitem{alessi2017}
R.~Alessi, J.-J. Marigo, C.~Maurini, and S.~Vidoli, ``Coupling damage and
  plasticity for a phase-field regularisation of brittle, cohesive and ductile
  fracture: One-dimensional examples,'' {\em International Journal of
  Mechanical Sciences}, 2017.
\newblock https://doi.org/10.1016/j.ijmecsci.2017.05.047.

\bibitem{shanthraj2016phase}
P.~Shanthraj, L.~Sharma, B.~Svendsen, F.~Roters, and D.~Raabe, ``A phase field
  model for damage in elasto-viscoplastic materials,'' {\em Computer Methods in
  Applied Mechanics and Engineering}, vol.~312, pp.~167--185, 2016.

\bibitem{choo18}
J.~Choo and W.~Sun, ``Coupled phase-field and plasticity modeling of geological
  materials: From brittle fracture to ductile flow,'' {\em Computer Methods in
  Applied Mechanics and Engineering}, vol.~330, pp.~1--32, 2018.

\bibitem{fang2019phase}
J.~Fang, C.~Wu, J.~Li, Q.~Liu, C.~Wu, G.~Sun, and Q.~Li, ``Phase field fracture
  in elasto-plastic solids: Variational formulation for multi-surface
  plasticity and effects of plastic yield surfaces and hardening,'' {\em
  International Journal of Mechanical Sciences}, vol.~156, pp.~382--396, 2019.

\bibitem{kruger2019porous}
M.~Kr{\"u}ger, M.~Dittmann, F.~Aldakheel, A.~H{\"a}rtel, P.~Wriggers, and
  C.~Hesch, ``Porous-ductile fracture in thermo-elasto-plastic solids with
  contact applications,'' {\em Computational Mechanics}, pp.~1--26, 2019.
\newblock https://doi.org/10.1007/s00466-019-01802-3.

\bibitem{nguyen2019multiscale}
L.~H. Nguyen and D.~Schillinger, ``The multiscale finite element method for
  nonlinear continuum localization problems at full fine-scale fidelity,
  illustrated through phase-field fracture and plasticity,'' {\em Journal of
  Computational Physics}, vol.~396, pp.~129--160, 2019.

\bibitem{dean2020phase}
A.~Dean, J.~Reinoso, N.~Jha, E.~Mahdi, and R.~Rolfes, ``A phase field approach
  for ductile fracture of short fibre reinforced composites,'' {\em Theoretical
  and Applied Fracture Mechanics}, p.~102495, 2020.
\newblock https://doi.org/10.1016/j.tafmec.2020.102495.

\bibitem{ali2021residual}
B.~Ali, Y.~Heider, and B.~Markert, ``Residual stresses in gas tungsten arc
  welding: a novel phase-field thermo-elastoplasticity modeling and parameter
  treatment framework,'' {\em Computational Mechanics}, pp.~1--23, 2021.

\bibitem{bryant2021phase}
E.~C. Bryant and W.~Sun, ``Phase field modeling of frictional slip with slip
  weakening/strengthening under non-isothermal conditions,'' {\em Computer
  Methods in Applied Mechanics and Engineering}, vol.~375, p.~113557, 2021.

\bibitem{storm2021comparative}
J.~Storm, M.~Pise, D.~Brands, J.~Schr{\"o}der, and M.~Kaliske, ``A comparative
  study of micro-mechanical models for fiber pullout behavior of reinforced
  high performance concrete,'' {\em Engineering Fracture Mechanics}, vol.~243,
  p.~107506, 2021.

\bibitem{ulloa2022micromechanics}
J.~Ulloa, J.~Wambacq, R.~Alessi, E.~Samaniego, G.~Degrande, and
  S.~Fran{\c{c}}ois, ``A micromechanics-based variational phase-field model for
  fracture in geomaterials with brittle-tensile and compressive-ductile
  behavior,'' {\em Journal of the Mechanics and Physics of Solids}, vol.~159,
  p.~104684, 2022.

\bibitem{khalil2022generalised}
Z.~Khalil, A.~Y. Elghazouli, and E.~Mart{\'\i}nez-Pa{\~n}eda, ``A generalised
  phase field model for fatigue crack growth in elastic--plastic solids with an
  efficient monolithic solver,'' {\em Computer Methods in Applied Mechanics and
  Engineering}, vol.~388, p.~114286, 2022.

\bibitem{guoliang1993monte}
J.~Guoliang, C.~Lin, and D.~Jiamei, ``{Monte Carlo finite element method of
  structure reliability analysis},'' {\em Reliability Engineering \& System
  Safety}, vol.~40, no.~1, pp.~77--83, 1993.

\bibitem{van1995use}
G.~Van~Vinckenroy and W.~De~Wilde, ``{The use of Monte Carlo techniques in
  statistical finite element methods for the determination of the structural
  behaviour of composite materials structural components},'' {\em Composite
  structures}, vol.~32, no.~1-4, pp.~247--253, 1995.

\bibitem{ganesh2021quasi}
M.~Ganesh, F.~Y. Kuo, and I.~H. Sloan, ``{Quasi-Monte Carlo finite element
  analysis for wave propagation in heterogeneous random media},'' {\em SIAM/ASA
  Journal on Uncertainty Quantification}, vol.~9, no.~1, pp.~106--134, 2021.

\bibitem{dick2013high}
J.~Dick, F.~Y. Kuo, and I.~H. Sloan, ``{High-dimensional integration: the
  quasi-Monte Carlo way},'' {\em Acta Numerica}, vol.~22, pp.~133--288, 2013.

\bibitem{kuo2012quasi}
F.~Y. Kuo, C.~Schwab, and I.~H. Sloan, ``{Quasi-Monte Carlo finite element
  methods for a class of elliptic partial differential equations with random
  coefficients},'' {\em SIAM Journal on Numerical Analysis}, vol.~50, no.~6,
  pp.~3351--3374, 2012.

\bibitem{graham2011quasi}
I.~G. Graham, F.~Y. Kuo, D.~Nuyens, R.~Scheichl, and I.~H. Sloan,
  ``{Quasi-Monte Carlo methods for elliptic PDEs with random coefficients and
  applications},'' {\em Journal of Computational Physics}, vol.~230, no.~10,
  pp.~3668--3694, 2011.

\bibitem{fang2022multilevel}
W.~Fang, Z.~Wang, M.~B. Giles, C.~H. Jackson, N.~J. Welton, C.~Andrieu, and
  H.~Thom, ``{Multilevel and Quasi Monte Carlo Methods for the Calculation of
  the Expected Value of Partial Perfect Information},'' {\em Medical Decision
  Making}, vol.~42, no.~2, pp.~168--181, 2022.

\bibitem{faustmann2021stability}
M.~Faustmann, J.~M. Melenk, and M.~Parvizi, ``{On the stability of Scott-Zhang
  type operators and application to multilevel preconditioning in fractional
  diffusion},'' {\em ESAIM: Mathematical Modelling and Numerical Analysis},
  vol.~55, no.~2, pp.~595--625, 2021.

\bibitem{kuo2017multilevel}
F.~Kuo, R.~Scheichl, C.~Schwab, I.~Sloan, and E.~Ullmann, ``{Multilevel
  quasi-Monte Carlo methods for lognormal diffusion problems},'' {\em
  Mathematics of Computation}, vol.~86, no.~308, pp.~2827--2860, 2017.

\bibitem{ben2020importance}
C.~Ben~Hammouda, N.~Ben~Rached, and R.~Tempone, ``{Importance sampling for a
  robust and efficient multilevel Monte Carlo estimator for stochastic reaction
  networks},'' {\em Statistics and Computing}, vol.~30, no.~6, pp.~1665--1689,
  2020.

\bibitem{kuo2015multi}
F.~Y. Kuo, C.~Schwab, and I.~H. Sloan, ``{Multi-level quasi-Monte Carlo finite
  element methods for a class of elliptic PDEs with random coefficients},''
  {\em Foundations of Computational Mathematics}, vol.~15, no.~2, pp.~411--449,
  2015.

\bibitem{chen2020insight}
Q.~Chen, H.~Zhu, H.~Li, J.~W. Ju, Z.~Jiang, and Z.~Yan, ``Insight into the
  inherent randomness of concrete properties using the stochastic
  micromechanics,'' {\em Probabilistic Engineering Mechanics}, vol.~61,
  p.~103064, 2020.

\bibitem{ibrahimbegovic2020reduced}
A.~Ibrahimbegovic, H.~G. Matthies, and E.~Karaveli{\'c}, ``Reduced model of
  macro-scale stochastic plasticity identification by bayesian inference:
  Application to quasi-brittle failure of concrete,'' {\em Computer Methods in
  Applied Mechanics and Engineering}, vol.~372, p.~113428, 2020.

\bibitem{gerasimov2020stochastic}
T.~Gerasimov, U.~R{\"o}mer, J.~Vond{\v{r}}ejc, H.~G. Matthies, and
  L.~De~Lorenzis, ``Stochastic phase-field modeling of brittle fracture:
  computing multiple crack patterns and their probabilities,'' {\em Computer
  Methods in Applied Mechanics and Engineering}, vol.~372, p.~113353, 2020.

\bibitem{ricoeur2022stochastic}
A.~Ricoeur, F.~Lindner, and K.~Zarjov, ``Stochastic aspects of crack deflection
  and crack path prediction in short fiber reinforced polymer matrix
  composites,'' {\em European Journal of Mechanics-A/Solids}, p.~104598, 2022.

\bibitem{colliat2007stochastic}
J.-B. Colliat, M.~Hautefeuille, A.~Ibrahimbegovic, and H.~G. Matthies,
  ``Stochastic approach to size effect in quasi-brittle materials,'' {\em
  Comptes Rendus M{\'e}canique}, vol.~335, no.~8, pp.~430--435, 2007.

\bibitem{dsouza2021non}
S.~M. Dsouza, T.~V. Mathew, I.~V. Singh, S.~Natarajan, {\em et~al.}, ``A
  non-intrusive stochastic phase field method for crack propagation in
  functionally graded materials,'' {\em Acta Mechanica}, vol.~232, no.~7,
  pp.~2555--2574, 2021.

\bibitem{matthies2005galerkin}
H.~G. Matthies and A.~Keese, ``Galerkin methods for linear and nonlinear
  elliptic stochastic partial differential equations,'' {\em Computer methods
  in applied mechanics and engineering}, vol.~194, no.~12-16, pp.~1295--1331,
  2005.

\bibitem{martin2022korali}
S.~M. Martin, D.~W{\"a}lchli, G.~Arampatzis, A.~E. Economides, P.~Karnakov, and
  P.~Koumoutsakos, ``Korali: Efficient and scalable software framework for
  bayesian uncertainty quantification and stochastic optimization,'' {\em
  Computer Methods in Applied Mechanics and Engineering}, vol.~389, p.~114264,
  2022.

\bibitem{pryse2017stochastic}
S.~Pryse and S.~Adhikari, ``Stochastic finite element response analysis using
  random eigenfunction expansion,'' {\em Computers \& Structures}, vol.~192,
  pp.~1--15, 2017.

\bibitem{matthies2008inelastic}
H.~G. Matthies and B.~V. Rosi{\'c}, ``Inelastic media under uncertainty:
  stochastic models and computational approaches,'' in {\em IUTAM Symposium on
  Theoretical, Computational and Modelling Aspects of Inelastic Media},
  pp.~185--194, Springer, 2008.

\bibitem{rosic2011stochastic}
B.~V. Rosic and H.~G. Matthies, ``Stochastic galerkin method for the
  elastoplasticity problem with uncertain parameters,'' in {\em Recent
  Developments and Innovative Applications in Computational Mechanics},
  pp.~303--310, Springer, 2011.

\bibitem{ben2016stochastic}
A.~Ben~Abdessalem, R.~Aza{\"\i}s, M.~Touzet-Cortina, A.~G{\'e}gout-Petit, and
  M.~Puiggali, ``Stochastic modelling and prediction of fatigue crack
  propagation using piecewise-deterministic markov processes,'' {\em
  Proceedings of the Institution of Mechanical Engineers, Part O: Journal of
  Risk and Reliability}, vol.~230, no.~4, pp.~405--416, 2016.

\bibitem{novak2005stochastic}
D.~Novak, M.~Vorechovsky, D.~Lehky, R.~Rusina, R.~Pukl, and V.~Cervenka,
  ``Stochastic nonlinear fracture mechanics finite element analysis of concrete
  structures,'' in {\em ICOSSAR05 Proc., 9th Int. Conf. on Structural Safety
  and Reliability}, pp.~781--788, Millpress, Rotterdam, The Netherlands, Rome,
  Italy, 2005.

\bibitem{lal2017stochastic}
A.~Lal, S.~B. Mulani, and R.~K. Kapania, ``Stochastic fracture response and
  crack growth analysis of laminated composite edge crack beams using extended
  finite element method,'' {\em International Journal of Applied Mechanics},
  vol.~9, no.~04, p.~1750061, 2017.

\bibitem{beck2013stochastic}
A.~T. Beck and W.~J. de~Santana~Gomes, ``Stochastic fracture mechanics using
  polynomial chaos,'' {\em Probabilistic Engineering Mechanics}, vol.~34,
  pp.~26--39, 2013.

\bibitem{junker2018analytical}
P.~Junker and J.~Nagel, ``An analytical approach to modeling the stochastic
  behavior of visco-elastic materials,'' {\em ZAMM-Journal of Applied
  Mathematics and Mechanics/Zeitschrift f{\"u}r Angewandte Mathematik und
  Mechanik}, vol.~98, no.~7, pp.~1249--1260, 2018.

\bibitem{junker2019relaxation}
P.~Junker and J.~Nagel, ``A relaxation approach to modeling the stochastic
  behavior of elastic materials,'' {\em European Journal of
  Mechanics-A/Solids}, vol.~73, pp.~192--203, 2019.

\bibitem{junker2020modeling}
P.~Junker and J.~Nagel, ``Modeling of viscoelastic structures with random
  material properties using time-separated stochastic mechanics,'' {\em
  International Journal for Numerical Methods in Engineering}, vol.~121, no.~2,
  pp.~308--333, 2020.

\bibitem{ghanem2003stochastic}
R.~G. Ghanem and P.~D. Spanos, {\em Stochastic finite elements: a spectral
  approach}.
\newblock Courier Corporation, 2003.

\bibitem{stefanou2009stochastic}
G.~Stefanou, ``The stochastic finite element method: past, present and
  future,'' {\em Computer methods in applied mechanics and engineering},
  vol.~198, no.~9-12, pp.~1031--1051, 2009.

\bibitem{pryse2021neumann}
S.~Pryse and S.~Adhikari, ``Neumann enriched polynomial chaos approach for
  stochastic finite element problems,'' {\em Probabilistic Engineering
  Mechanics}, vol.~66, p.~103157, 2021.

\bibitem{reddy2008stochastic}
R.~Reddy and B.~Rao, ``Stochastic fracture mechanics by fractal finite element
  method,'' {\em Computer Methods in Applied Mechanics and Engineering},
  vol.~198, no.~3-4, pp.~459--474, 2008.

\bibitem{khodadadian2020bayesian}
A.~Khodadadian, N.~Noii, M.~Parvizi, M.~Abbaszadeh, T.~Wick, and C.~Heitzinger,
  ``A bayesian estimation method for variational phase-field fracture
  problems,'' {\em Computational Mechanics}, vol.~66, pp.~827--849, 2020.

\bibitem{noii2021bayesian}
N.~Noii, A.~Khodadadian, J.~Ulloa, F.~Aldakheel, T.~Wick, S.~Fran{\c{c}}ois,
  and P.~Wriggers, ``Bayesian inversion for unified ductile phase-field
  fracture,'' {\em Computational Mechanics}, vol.~68, no.~4, pp.~943--980,
  2021.

\bibitem{abbaszadeh2021reduced}
M.~Abbaszadeh, M.~Dehghan, A.~Khodadadian, N.~Noii, C.~Heitzinger, and T.~Wick,
  ``A reduced-order variational multiscale interpolating element free galerkin
  technique based on proper orthogonal decomposition for solving navier--stokes
  equations coupled with a heat transfer equation: Nonstationary incompressible
  boussinesq equations,'' {\em Journal of Computational Physics}, vol.~426,
  p.~109875, 2021.

\bibitem{noii2022bayesian}
N.~Noii, A.~Khodadadian, and T.~Wick, ``Bayesian inversion using global-local
  forward models applied to fracture propagation in porous media,'' {\em
  International Journal for Multiscale Computational Engineering}, 2022.

\bibitem{ALDAKHEEL2021114175}
F.~Aldakheel, N.~Noii, T.~Wick, O.~Allix, and P.~Wriggers, ``Multilevel
  global-local techniques for adaptive ductile phase-field fracture,'' {\em
  Computer Methods in Applied Mechanics and Engineering}, vol.~387, p.~114175,
  2021.

\bibitem{hughes2012finite}
T.~J. Hughes, {\em The finite element method: linear static and dynamic finite
  element analysis}.
\newblock Courier Corporation, 2012.

\bibitem{taghizadeh2017optimal}
L.~Taghizadeh, A.~Khodadadian, and C.~Heitzinger, ``{The optimal multilevel
  Monte-Carlo approximation of the stochastic drift--diffusion-Poisson
  system},'' {\em Computer Methods in Applied Mechanics and Engineering},
  vol.~318, pp.~739--761, 2017.

\bibitem{khodadadian2020adaptive}
A.~Khodadadian, M.~Parvizi, and C.~Heitzinger, ``{An adaptive multilevel Monte
  Carlo algorithm for the stochastic drift--diffusion--Poisson system},'' {\em
  Computer Methods in Applied Mechanics and Engineering}, vol.~368, p.~113163,
  2020.

\bibitem{miehe2011}
C.~Miehe, ``A multi-field incremental variational framework for
  gradient-extended standard dissipative solids,'' {\em Journal of the
  Mechanics and Physics of Solids}, vol.~59, no.~4, pp.~898--923, 2011.

\bibitem{wittmann1985simulation}
F.~Wittmann, P.~Roelfstra, and H.~Sadouki, ``Simulation and analysis of
  composite structures,'' {\em Materials Science and Engineering}, vol.~68,
  no.~2, pp.~239--248, 1985.

\bibitem{wang1999mesoscopic}
Z.~Wang, A.~Kwan, and H.~Chan, ``Mesoscopic study of concrete i: generation of
  random aggregate structure and finite element mesh,'' {\em Computers \&
  Structures}, vol.~70, no.~5, pp.~533--544, 1999.

\bibitem{wriggers2006mesoscale}
P.~Wriggers and S.~Moftah, ``Mesoscale models for concrete: Homogenisation and
  damage behaviour,'' {\em Finite Elements in Analysis and Design}, vol.~42,
  no.~7, pp.~623--636, 2006.

\end{thebibliography}
	 	\end{spacing}}
	
\end{document}